\documentclass[11pt]{article}
\usepackage{amsmath,amsfonts,amssymb,verbatim,latexsym,cases,mathdots,colordvi,url,extarrows}
\usepackage[linesnumbered,ruled,vlined]{algorithm2e}
\usepackage{relsize}
\usepackage{tikz}
\usepackage{graphicx}
\usepackage{epstopdf}
\usepackage{mathrsfs}
\usepackage{color}
\usepackage{amsthm}
\usepackage{enumerate}
\usepackage{indentfirst}
\usepackage{float}

\renewcommand{\Box}{\rule{2.2mm}{2.2mm}}

\def\beginproof{\par\noindent {\bf Proof.}\ \ }
\def\endproof{\hskip .5cm $\Box$ \vskip .5cm}

\begin{document}
\newtheorem{theorem}{Theorem}[section]
\newtheorem{definition}[theorem]{Definition}
\newtheorem{lemma}[theorem]{Lemma}
\newtheorem{assump}[theorem]{Assumption}
\newtheorem{proposition}[theorem]{Proposition}
\newtheorem{corollary}[theorem]{Corollary}
\newtheorem{example}[theorem]{Example}
\newtheorem{remark}[theorem]{Remark}
\newtheorem{refer}{Reference}
\begin{sloppypar}
\title{New second-order optimality conditions for  directional optimality of a general set-constrained optimization problem\thanks{The alphabetical order of the paper indicates the equal contribution to the paper.}}
\author{Wei Ouyang\thanks{School of Mathematics, Yunnan Normal University, Kunming 650500, People's Republic of China; Yunnan Key Laboratory of Modern Analytical Mathematics and Applications, Kunming 650500, China. The research of this author was partially supported by the National Natural Science Foundation of the People's Republic of China [Grant 12261109] and the Basic Research Program of Yunnan Province
[Grant 202301AT070080]. Email: weiouyangxe@hotmail.com. } \and Jane J. Ye \thanks{Corresponding author. Department of Mathematics and Statistics, University of Victoria, Canada. The research of this author was partially
supported by NSERC. Email: janeye@uvic.ca.} \and Binbin Zhang\thanks{School of Science, Kunming University of Science and Technology, Kunming 650500, People's Republic of China.  Email: bbzhang@kust.edu.cn.}}
  \date{}
\maketitle 
\begin{abstract} 
	In this paper we derive new second-order optimality conditions for a very general set-constrained optimization problem where the underlying set may be nononvex. We consider local optimality in specific directions  in pursuit of developing these new  optimality conditions. Utilizing the classical and/or the lower generalized support function, we obtain new second-order necessary and sufficient conditions for local optimality of general nonconvex constrained optimization problems in given directions via both the corresponding asymptotic second-order tangent cone and outer second-order tangent set. Our results do not require  convexity and/or nonemptiness of the outer second-order tangent set. This is an important improvement to other results in the literature since the outer second-order tangent set can be nonconvex and  empty even when the set is convex.

\vskip 10 true pt
\noindent {\bf Key words.}\quad first-order and second-order optimality condition, directional limiting normal cone, outer second-order tangent set, asymptotic second-order tangent cone, support function, lower generalized support function
\vskip 10 true pt


\end{abstract}

\newpage
\section{Introduction}
In this paper, we consider the general set-constrained optimization problem in the form
\begin{equation}\label{op}
\min f(x)\,\,\, {\rm s.t.}\,\,\, g(x)\in K,
\end{equation}
where $f:\mathbb{R}^n\rightarrow \mathbb{R}$, $g:\mathbb{R}^n\rightarrow \mathbb{R}^m$ and  $K$ is a closed subset of $\mathbb{R}^m$. Throughout the paper we assume that $f,g$ are twice continuously differentiable  unless otherwise specified. For convenience, we denote by $C:=g^{-1}(K)$ the feasible region of problem (\ref{op}), $T_K(\cdot)$ the tangent cone to $K$ and {$C(\bar x):=\{d\in \mathbb{R}^n \mid \nabla g(\bar x)d \in T_K(g(\bar x)), \nabla f(\bar x)d\leq 0\}$} the critical cone at a feasible solution $\bar x$. In fact,  the results of this paper can be suitably modified to include the more general case where the image space  is  a finite dimensional Hilbert space as in Guo et al. \cite{guo2013mathematical}.

When $K$ is convex, problem (\ref{op}) includes the classical optimization problem classes  such as the nonlinear program, the nonlinear conic  program (when $K$ is a convex cone),  the second-order cone program and the semidefinite cone program as special cases. The possibility of $K$ being nonconvex  in \eqref{op} extends the problem (\ref{op}) to include much more general problems and provides a generic framework for a number of important optimization problem classes such as the mathematical program with equilibrium constraints \cite{lzq96}, the mathematical program with second-order cone complementarity constraints (Ye and Zhou \cite{yjj16,yjj18}) and the mathematical program with semidefinite cone complementarity constraints (Ding et al. \cite{dc14}). 

In the case of a nonlinear program (e.g. when $K= \mathbb{R}_-^{m}$), recall that if $\bar x$ is a local optimal solution of  problem \eqref{op}, then under certain constraint qualification (e.g. Mangasarian-Fromovitz constraint qualification), there is a multiplier $\lambda$ fullfilling the first order necessary optimality condition such that  the second order necessary optimality condition holds:
\begin{equation}\label{spcNLP}
\nabla^2_{xx}L(\bar x,\lambda)(d,d)\geq 0\qquad  \forall d\in C(\bar x),
\end{equation}
 where $L(x,\lambda):=f(x)+\langle\lambda,g(x)\rangle$ denotes the Lagrangian. {Conversely}, if (\ref{spcNLP}) holds with strict inequality ``$>$'' for all nonzero critical direction $d\in C(\bar x)\setminus \{0\}$, for certain feasible solution $\bar x$, a certain multiplier $\lambda$,  then $\bar x$ must be a locally optimal solution. This kind of second-order necessary and sufficient optimality condition is  of ``no gap'' meaning that the only difference between the sufficent and the necessary optimality condition is that 
the inequality ``$\geq $'' in the  necessary optimality condition is changed to the strict inequality in the sufficient condition.

Recently, second-order variational analysis and optimality conditions for problem in the form \eqref{op} have been developed rapidly; see \cite{ben23a,ben23,fran24,gfer19,gfer16,ye22} and the references therein.
For problem \eqref{op} with $K$ being convex, there is complete theory of optimality conditions in the existing literature (see Bonnans et al. \cite{bon99} and the monograph of Bonnans and Shapiro \cite{shap}). The second-order necessary  optimality conditions are derived within the framework of convex analysis with the following form (cf. \cite[Theorem 3.45]{shap}): If the Robinson constraint qualification (Robinson CQ) holds at a local minimizer $\bar x$, then for every critical direction $d\in C(\bar x)$ and every convex subset $K(d)$ of the second-order tangent set $T_K^2(g(\bar x);\nabla g(\bar x)d)$, there is a multiplier $\lambda$ fulfilling the first-order necessary optimality condition such that 
\begin{equation}\label{spc}
\nabla^2_{xx}L(\bar x,\lambda)(d,d)-\sigma_{K(d)}(\lambda)\geq 0,
\end{equation}
where  $\sigma_{K(d)}(\cdot)$ represents the classic support function to set ${K(d)}$.
In particular if $K$ and  $T_K^2(g(\bar x);\nabla g(\bar x)d)$ are both convex, then $K(d)$ in  (\ref{spc}) can be taken as $T_K^2(g(\bar x);\nabla g(\bar x)d)$ and (\ref{spc}) becomes the following condition
\begin{equation}\label{spcnew}
\nabla^2_{xx}L(\bar x,\lambda)(d,d)-\sigma_{T_K^2(g(\bar x);\nabla g(\bar x)d)}(\lambda)\geq 0.
\end{equation}
 Unfortunately, however,  the second-order tangent set of a convex set is generally not convex and hence even when $K$ is convex, (\ref{spcnew}) may not hold as a second-order necessary optimality condition.

For the general case where $K$ may be  nonconvex, Gfrerer et al. \cite[Theorem 2]{ye22} has obtained the following second-order necessary optimality condition:  Let $\bar x$ be a local minimum. If the directional metric subregularity constraint qualification (MSCQ), which is weaker than  Robinson CQ, holds at $\bar x$ in direction $d$ with 
$d\in C(\bar x)$ and $T_K^2(g(\bar x);\nabla g(\bar x)d)\not =\emptyset$, then there exists a Mordukhovich (M-) multiplier $\lambda$ in direction $d$ such that
\begin{equation}\label{spc1}
\nabla^2_{xx}L(\bar x,\lambda)(d,d)-\hat{\sigma}_{T_K^2(g(\bar x);\nabla g(\bar x)d)}(\lambda)\geq 0,
\end{equation}
where $\hat\sigma_C(\cdot)$  is the so-called lower generalized support function (see Definition \ref{gsf}).  The function $\hat\sigma_C$ is indeed an extension of the classic support function and it is ``lower'' than the support function since it  has property that 
$\hat\sigma_C\leq \sigma_C$, and the equality holds when $C$ is closed and conex. Hence (\ref{spc1}) is a weaker condition than (\ref{spcnew}).  In general if $K$ is nonconvex we can only obtain the weaker condition (\ref{spc1}). However it was shown in \cite[Proposition 3.46]{shap} and Gfrerer et al. \cite[Corollary 5]{ye22} that the stronger condition (\ref{spcnew}) still holds even when $K$ is nonconvex if the multipliers considered are unique.
Under the assumption that the set $K$ is  outer second-order regular and the second-order tangent set
 $T_K^2(g(\bar x);\nabla g(\bar x)d)$ is convex, by the classical result    \cite[Theorem 3.86]{shap} and \cite[Theorem 4]{ye22}, the second order necessary optimality condition (\ref{spcnew}) becomes the following ``no gap''  sufficient optimality condition: Let  $\bar x$  be a feasible solution of problem \eqref{op}. If for each nonzero critical direction $d$ there is $ \lambda$ fullfiling the first order necessary optimality condition such that
\begin{equation}\label{spcsuffcl}  \nabla_{xx}^2L(\bar x,\lambda)(d,d)-\sigma_{T^2_K (g(\bar x);\nabla g({\bar x})d)}(\lambda)>0,
\end{equation}
then  $\bar x$ is a local optimal solution of problem \eqref{op}. Recently sufficient optimality conditions with the second-order subderivative of the indicator function $\delta_K$ are developed first in Mohammadi et al. \cite[Proposition 7.3]{moha}, under the convexity and   the parabolical derivability of set $K$ and with these assumptions  removed in \cite[Theorem 3.3]{ben23a} from which we obtain a second-order  sufficient optimality condition as follows: Let $\bar x$ be a feasible solution. Suppose that for every nonzero critical direction $d$ there exists $ \lambda$ fullfiling the first order necessary optimality condition such that 
\begin{equation}\label{spcsuff}
\nabla^2_{xx} L(\bar x,\lambda)(d,d)+d^2\delta_K(g(\bar x);\lambda)(\nabla g(\bar x)d)> 0.
\end{equation} Then $\bar x$ is a local minimizer.

From our discussions so far, we see that the second-order optimality conditions rely on the three second-order objects $\hat\sigma_{T^2_K(g(\bar x);\nabla g(\bar x)d))}(\lambda)$, $\sigma_{T^2_K(g(\bar x);\nabla g(\bar x)d))}(\lambda)$ and ${\rm d}^2\delta_{K}(g(\bar x);\lambda)(\nabla g(\bar x)d))$, each of them describing in some way the curvature of the set $K$. It is known by \cite[Proposition 2.18]{ben23} that these three second-order objectives are linked together by the inequalities
\begin{equation}\label{3_curv_terms_rel}
{{\rm d}^2\delta_K(g(\bar x);\lambda)(\nabla g(\bar x)d)\leq  -\sigma_{T^2_K(g(\bar x);\nabla g(\bar x)d)}(\lambda)\leq -\hat{\sigma}_{T^2_K(g(\bar x);\nabla g(\bar x)d)}(\lambda).}
\end{equation}
Since the equalities in general do not hold, (\ref{spcsuff}) cannot be  a ``no gap'' sufficient condition of the necessary condition (\ref{spc1}).

For an empty set, by convention $\sigma_\emptyset(\lambda):=-\infty$ and by definition $\hat \sigma_\emptyset(\lambda):=-\infty$ as well. Hence when 
$T_K^2(g(\bar x);\nabla g(\bar x)d)=\emptyset$, conditions 
 (\ref{spc}), (\ref{spcnew}) and  (\ref{spc1}) hold trivially. This means that if $T_K^2(g(\bar x);\nabla g(\bar x)d)=\emptyset$, conditions (\ref{spc}),
(\ref{spcnew}) and  (\ref{spc1}) do not provide any useful information. Moreover the outer second-order regularity of $K$ ensures that nonemptiness of the set  $T_K^2(g(\bar x);\nabla g(\bar x)d)$ (see \cite[Page 202]{shap}).

In this paper we aim at deriving second-order optimality conditions which do not require convexity of set $K$, the  outer second-order regularity of set $K$ and the nonemptiness of the  second-order tangent set
$T_K^2(g(\bar x);\nabla g(\bar x)d)$. Moreover we wish to derive second-order necessary conditions { that are as close} to  the corresponding second-order sufficient condition as possible. In \cite{pen98}, Penot  introduced the notion of asymptotic second-order tangent cone $T_K^{''}(g(\bar x);\nabla g(\bar x)d)$ and used it in the study of second-order optimality conditions for  optimization problems with an abstract set constraint (i.e., the case when the constraint is $x\in K$). He observed the following important fact. The set of nonzero elements in asymptotic second-order tangent cone $T_K^{''}(g(\bar x);\nabla g(\bar x)d)$ and the second-order tangent set $T_K^2(g(\bar x);\nabla g(\bar x)d)$ cannot be empty simultaneously and had used both these sets in the second-order optimality conditions. Taking further Penot's idea, in this paper we obtain second-order optimality conditions using both the second-order tangent set and the asymptotic second-order tangent cone and hence no longer require the nonemptiness of the  set
$T_K^2(g(\bar x);\nabla g(\bar x)d)$ in our optimality conditions.

We observe that the following relationships always hold:
\begin{eqnarray}
\Omega &:=& \nabla g(\bar x)(T_C^2(\bar x;d))+\nabla^2 g(\bar x)(d,d)\subseteq T_K^2(g(\bar x);\nabla g(\bar x)d),\label{strin}\\
\Theta  &:=& \nabla g(\bar x)(T_{C}^{''}(\bar x;d))\subseteq T_K^{''}(g(\bar x);\nabla g(\bar x)d),\label{strintheta}
\end{eqnarray}
and the strict inclusions may hold (see Proposition \ref{qc} and Example \ref{fg3a}).
We succeed in showing that  (see {Theorem \ref{nss}), if problem \eqref{op} attains its local minimum  at $\bar x$ 
{and there is no descent direction, then  for any $\lambda$ with $\nabla_x L(\bar x,\lambda)=0$ and any critical direction $d$,}  the following conditions hold:
\begin{equation}\label{spc3theta}
\aligned
& \nabla^2_{xx}L(\bar x,\lambda)(d,d)-\sigma_{\Omega}(\lambda)\geq 0
 \mbox{ and }  \sigma_{\Theta}(\lambda)\leq  0.
 \endaligned
\end{equation}
Moreover by Theorem \ref{sff},  we can obtain the following ``no gap'' sufficient optimality condition of the necessary optimality condition (\ref{spc3theta}). Let $\bar x$ be a feasible solution.   Suppose there is $\lambda$ with $\nabla_xL(\bar x,\lambda)=0$ and there is no descent direction. Suppose for all nonzero critical directions $d$ the following conditions hold:
\begin{equation}\label{spc3s}
\nabla^2_{xx}L(\bar x,\lambda)(d,d)-\sigma_\Omega(\lambda)>0
\ \ \text{and}\ \ \sigma_{\Theta \backslash \{0\}}(\lambda)< 0.
\end{equation}
Then $\bar x$ is a local optimal solution. 

Since the objects $\Omega$ and $\Theta$ are implicit,  conditions (\ref{spc3theta}) and (\ref{spc3s}) themselves are of limited practical use.
However, by Proposition \ref{remark4.6}, if $\bar x$ is a local minimum of problem \eqref{op}  and  the directional MSCQ holds at $\bar x$ in a  critical direction $d$, then we have the following relationships:
\begin{itemize}
\item[(i)] 
There exists an M-multiplier $\lambda$ in direction $d$ such that 
$\hat{\sigma}_{{T^{''}_K (g(\bar x);\nabla g(\bar x)d)}}(\lambda)\leq 
\sigma_{{\Theta}}(\lambda).
$
\item[(ii)] If $T^2_K (g(\bar x);\nabla g(\bar x)d)\not =\emptyset$, then there exists an M-multiplier $\lambda$ in direction $d$ such that 
$\hat{\sigma}_{{T^2_K (g(\bar x);\nabla g(\bar x)d)}}(\lambda)\leq \sigma_\Omega(\lambda).$
\end{itemize}
These relationships together  with conditions (\ref{spc3theta}) 
can  be used to derive second-order necessary optimality condition involving $\hat\sigma_{T^2_K (g(\bar x);\nabla g(\bar x)d)}(\lambda)$.
Indeed  we have improved Gfrerer et al. \cite[Theorem 2]{ye22} to the following second-order necessary optimality condition (Theorem  \ref{cclmby} with $A=\mathbb{R}^n$ and $B=\mathbb{R}^n$): Let $\bar x$ be a local minimum and $d\in C(\bar x)$ be a critical direction. If  the directional MSCQ holds at $\bar x$ in direction $d$, then the following conditions hold:
\begin{itemize}
\item[(i)] 
There exists an M-multiplier $\lambda$ in direction $d$ such that
$$\hat{\sigma}_{{T^{''}_K (g(\bar x);\nabla g(\bar x)d)}}(\lambda)\leq 0.$$
\item[(ii)] If $T^2_K (g(\bar x);\nabla g(\bar x)d)\not =\emptyset$, then there exists an M-multiplier $\lambda$ in direction $d$ such that 
 $$\nabla^2_{xx}L(\bar x,\lambda)(d,d)-\hat{\sigma}_{{T^2_K (g(\bar x);\nabla g(\bar x)d)}}(\lambda)\geq 0.$$ 
\end{itemize}
Furthermore in Corollary \ref{fffd}, we have provided an improvement of the classical second-order necessary optimality condition (\ref{spc}) in that $K$ can be nonconvex and $T_K^2(g(\bar x);\nabla g(\bar x)d)$ can be empty, and  by  Corollary \ref{fffdss}, under the directional nondegeneracy condition which is stronger than the directional MSCQ, the lower generalized support function $\hat{\sigma}$ in the above conditions can be replaced by the support function $\sigma$.
Similarly by Corollary \ref{cor4.15}, we have improved the classical second-order sufficent optimality condition (\ref{spcsuffcl})
to the following second-order sufficient optimality condition: Let  $\bar x$  be a feasible solution of problem \eqref{op} and assume that there is no descent direction from $\bar x$.
Assume that for each nonzero critical direction  there is $\lambda\in \mathbb{R}^m$ with $\nabla_x L(\bar x,\lambda)=0$, such that
\begin{enumerate}[{\rm (i)}]
\item $\langle\lambda,v\rangle< 0,\quad\forall v\in  T''_K (g(\bar x);\nabla g({\bar x})d)  \cap\nabla g(\bar x) \left (\{d\}^\bot\backslash\{0\}\right )$;
\item $ \nabla_{xx}^2L(\bar x,\lambda)(d,d)-\sigma_{T^2_K (g(\bar x);\nabla g({\bar x})d)}(\lambda)>0$.
\end{enumerate}
Then, $\bar x$ is a local optimal solution of problem \eqref{op}. 
Note that unlike in the classical result \cite[Theorem 3.86]{shap} and \cite[Theorem 4]{ye22}, in our new second-order sufficient optimality condition, $K$ is not required to be outer second-order regular and $T^2_K (g(\bar x);\nabla g(\bar x)d)$ is not required to be convex.

In practical applications, it may be useful to study  local minima or maxima in certain directions.
For example in the exact line search of a gradient method, to find a step size, one may minimize the objective function in a certain  descent direction. Sometimes, information on a local minimum along a certain direction may have information related to the non-directional local minimum. For example, if we know a point is not a local minimum along a non-descent direction then it cannot be a local minimum in the classic sense as well.  To this end, we introduce in this paper the notion of  local minima in a certain direction. For a given reference point $\bar x$, if a local minimum value is obtained in each critical direction, then it reduces to the notion of local minimum in the classic sense.

Unlike the directional necessary and/or sufficient conditions obtained for local optimality (along all directions of a reference point) in the existing literature, in this paper, we  develop new optimality conditions for problem \eqref{op} to obtain a local minimum value in a certain direction.

The main contributions of this paper are as follows.
\begin{itemize}
\item We introduce the concept of a  local optimal solution in a given direction and derive first-order and second-order optimality conditions for a {directional} local optimum. This provides a more flexible framework since one can derive the corresponding first-order and second-order optimality conditions for a standard local optimum from a directional local optimum if the given directions are the set of all critical directions.
\item Unlike classical  second-order optimality conditions, our new second-order optimality condition for local optimality of nonconvex problem \eqref{op} in given directions  exhibits the exact correspondence between the direction of local optimality and the direction in optimality conditions.
\item We obtain second-order optimality conditions without additional restrictions (such as convexity of set $K$, outer second-order regularity on the set $K$, nonemptiness and/or of the  second-order tangent set
$T_K^2(g(\bar x);\nabla g(\bar x)d)$ etc.). Moreover our new second-order sufficient optimality condition (\ref{spc3s}) is  of ``no gap'' to its corresponding second-order necessary optimality condition (\ref{spc3theta}).
\end{itemize}

We organize our paper as follows. Section 2 contains the basic notations and preliminary results.
In Section 3, we derive new first- and second-order necessary optimality conditions for nonconvex optimization problem with an abstract set constraint in given directions. In Section 4, we develop second-order optimality conditions for general nonconvex optimization problem \eqref{op} in given directions.

\section{Notations and Preliminary Results}
\label{sec:notation}
In this section, we provide the basic notations and fundamental facts in variational analysis which are used throughout the paper and develop some preliminary results.

We denote the closed unit ball of $\mathbb{R}^n$ by $B_{\mathbb{R}^n}$.
$B_r(x)$ represents the closed ball with center $x\in \mathbb{R}^n$ and radius $r>0$. For a set $S\subseteq \mathbb{R}^n$, we denote its interior, closure, boundary, convex hull and {linear hull} by ${\rm int}(S), {\rm cl}(S), {\rm bd}(S)$, ${\rm co}(S)$ and ${\rm span}(S)$ respectively. Let $S^\circ$ and $\sigma_S(x)$ stand for the polar cone
and the support function of $S$, respectively, that is, $S^\circ:=\{v\in \mathbb{R}^n \mid \langle v,x\rangle\leq 0\ \text{for\  all}\  x\in S\}$ and $\sigma_S(x):=\sup_{u\in S}\langle x,u\rangle$ for all $x\in \mathbb{R}^n$. {Let ${\rm cone}(S)$ denote the cone generated by $S$.}
Let $d(x,S):=\inf_{u\in S}\|x-u\|$ denote the point-to-set distance from $x$ to $S$ (in the usual convention, the infimum of the empty set equals $+\infty$).  Let $o: \mathbb{R}_+\rightarrow \mathbb{R}^n$ denote {a} mapping satisfying $o(t)/t\rightarrow 0$ as $t\downarrow 0$. For $u\in \mathbb{R}^n$, denote by $\{u\}^\bot$ the orthogonal complement of the linear space generated by $u$. For a vector-valued mapping $h:\mathbb{R}^n\rightarrow\mathbb{R}^m$ and sets $A\subseteq \mathbb{R}^n$, $B\subseteq \mathbb{R}^m$, we denote $h(A):=\{h(x)|x\in A\}$ and $h^{-1}(B):=\{x| h(x)\in B\}$.

For a vector-valued mapping $g: \mathbb{R}^n\rightarrow \mathbb{R}^m$, we denote by $\nabla g(x)\in \mathbb{R}^{m\times n}$ the Jacobian of $g$ at $x$. {When $g$ is a scalar mapping (i.e. $m=1$), we denote by $\nabla g(x)$ as either the gradient vector of $g$ or the Jacobian of $g$ depending on the context} and $\nabla^2g(x)$ is the Hessian at $x$. 

For the mapping $g: \mathbb{R}^n\rightarrow \mathbb{R}^m$, the second-order derivative of $g$ at $x$ is denoted by $\nabla^2 g(x)$ and is defined as follows:
$$u^T\nabla^2 g(x):=\lim_{t\rightarrow 0}\frac{\nabla g(x+tu)-\nabla g(x)}{t},\quad\forall u\in \mathbb{R}^n.$$
Hence, we have
$$\nabla^2 g(x)(u,u):=u^T\nabla^2 g(x)u=(u^T\nabla^2 g_1(x)u,\ldots,u^T\nabla^2 g_m(x)u),\quad\forall u\in \mathbb{R}^n.$$

For a set-valued mapping $M:\mathbb{R}^n\rightrightarrows \mathbb{R}^m$, its graph is defined by
${\rm gph} M:=\{(u,v)\in \mathbb{R}^n\times \mathbb{R}^m \mid v\in M(u)\}$. The inverse mapping $M^{-1}:\mathbb{R}^m\rightrightarrows \mathbb{R}^n$ is defined by $M^{-1}(v)=\{u\in \mathbb{R}^n \mid v\in M(u)\}$ for all $v\in \mathbb{R}^m$.
We denote by $\limsup_{u'\rightarrow u}M(u')$ {and $\liminf_{u'\rightarrow u}M(u')$} the
Painlev\'{e}-Kuratowski upper {and lower }limit respectively, that is,
$$\limsup\limits_{u'\rightarrow u}M(u'):=\{v\in \mathbb{R}^m \mid \exists u_k\rightarrow u,v_k\rightarrow v\,\,{\rm such\, that}\,\,v_k\in M(u_k)\},$$
\begin{eqnarray*}
\lefteqn{\liminf\limits_{u'\rightarrow u}M(u'):=}\\
&&\{v\in \mathbb{R}^m \mid \forall u_k\rightarrow u, \exists v_k\rightarrow v\,\,{\rm such\, that}\,\,v_k\in M(u_k) \mbox{ for sufficiently large  } k\}.
\end{eqnarray*}

\subsection{Tangent and Normal Cones}
For a closed subset $S$ of $\mathbb{R}^n$ and a point $\bar x \in S$, the regular/Clarke and Bouligand-Severi tangent/contingent cone to $S$ at $\bar x$ is defined, respectively, by
$$\hat T_S(\bar x):=\liminf\limits_{t\downarrow 0,x\xrightarrow[]{S}\bar x}\frac{S-x}{t}=\{d\in \mathbb{R}^n \mid \forall t_k\downarrow 0, x_k\xrightarrow[]{S}\bar x,\exists d_k\rightarrow d \,\,{\rm with}\,\, x_k+t_kd_k\in S\}$$
and
$$T_S(\bar x):=\limsup\limits_{t\downarrow 0}\frac{S-\bar x}{t}=\{d\in \mathbb{R}^n \mid \exists t_k\downarrow 0, d_k\rightarrow d \,\,{\rm with}\,\, \bar x+t_kd_k\in S\},$$
where $x\xrightarrow[]{S}\bar x$ represents the convergence of $x$ to $\bar x$ with $x\in S$.
The Fr\'{e}chet normal cone of $S$ to $\bar x$ is defined by
$$\hat N_S(\bar x):=\left\{v\in \mathbb{R}^n \mid \limsup_{x\xrightarrow[]{S}\bar x}\frac{\langle v,x-\bar x\rangle}{\|x-\bar x\|}\leq 0\right\}.$$
It is {well known} that $\hat N_S(\bar x)=(T_S(\bar x))^\circ$.

Let $N_S(\bar x)$ denote the limiting/Mordukhovich/basic normal cone of $S$ at $\bar x$, that is,
$$N_S(\bar x):=\limsup_{x\xrightarrow[]{S}\bar x}{ \hat{N}_S(x)}.$$ {By convention for any $\bar x\not\in S$, $N_S(\bar x)$ is defined as the empty set.}
The limiting normal cone is, in general, nonconvex, whereas the Fr\'{e}chet normal cone is always convex (see,
e.g., Mordukhovich {\cite[Definition 1.1]{boris}} and Rockafellar and Wets {\cite[Proposition 6.5]{rock}}). In the case of $S$ being convex, the Fr\'{e}chet and limiting normal cones coincide with the normal cone in the sense of convex analysis, that is,
$$N_S(\bar x):=\{v\in \mathbb{R}^n \mid \langle v,x-\bar x\rangle\leq 0\ \text{for\  all}\  x\in S\}.$$

The directional limiting normal cone was introduced in Ginchev and Mordukhovich \cite{ginc} and was extended to general Banach spaces in Gfrerer \cite{cfe13}. For a given direction $d\in \mathbb{R}^n$, the limiting normal cone to $S$ in direction $d$ at $\bar x$ is defined by
\begin{eqnarray*}
{N_S(\bar x;d):=\limsup_{t\downarrow 0,d'\rightarrow d}\hat{N}_S(\bar x+td').}
\end{eqnarray*}
From the definition, it is easy to observe that $N_S(\bar x;d)=\emptyset$ if $d\not\in T_S(\bar x)$, $N_S(\bar x;d)\subseteq N_S(\bar x)$ and $N_S(\bar x;0)= N_S(\bar x)$. When $S$ is convex and $d\in T_S(\bar x)$, we have the following relationship for directional normal cone:
\begin{equation}\label{ccc}
N_S(\bar x;d)=N_S(\bar x)\cap\{d\}^\bot=N_{T_S(\bar x)}(d).
\end{equation}
For more details, please refer to Gfrerer \cite[Lemma 2.1]{gfer14}.

Recently, the directional regular/Clarke tangent cone and directional Clarke normal cone were introduced by Gfrerer et al. \cite{ye22}. Let $\hat T_S(\bar x;d)$ denote the directional regular/Clarke tangent cone of $S$ at $\bar x$ in direction $d$, that is,
\begin{equation*}
\aligned
\hat T_S(\bar x;d)&:=\liminf_{t\downarrow 0,d'\rightarrow d, \bar x+td'\in S}T_S(\bar x+td').
\\
\endaligned
\end{equation*}
Let $N^c_S(\bar x;d):={\rm cl}\ {\rm co}\ (N_S(\bar x;d))$ denote the directional Clarke normal cone of $S$ at $\bar x$ in direction $d$.
It is easy to see from the definition that, for a closed set S, the directional version of the regular tangent cone
contains the nondirectional one, and it coincides with the nondirectional one when the direction is equal to zero,
that is, $\hat T_S(\bar x)\subset\hat T_S(\bar x;d)$ and $\hat T_S(\bar x)=\hat T_S(\bar x;0)$. It is shown in \cite[Proposition 3]{ye22} that, $\hat T_S(\bar x;d)$  is closed and convex and $(\hat T_S(\bar x;d))^\circ=N^c_S(\bar x;d)$.

For a proper lower semi-continuous function $\varphi:\mathbb{R}^n\rightarrow \mathbb{R}\cup \{+\infty\}$, let ${\rm dom }(\varphi)$ and ${\rm epi}(\varphi)$ denote the domain and the epigraph of $\varphi$, respectively, that is,
$${\rm dom} (\varphi):=\{x\in \mathbb{R}^n \mid \,\,\varphi(x)<+\infty\}\,\,{\rm and}\,\,{\rm epi}(\varphi):=\{(x,r)\in \mathbb{R}^n\times \mathbb{R} \mid \varphi(x)\leq r\}.$$
Let $\hat \partial\varphi(x)$ and $\partial\varphi(x)$ denote the Fr\'{e}chet and the limiting/Mordukhovich subdifferential of $\varphi$ at $x$, respectively, that is
$$\hat \partial\varphi(x):=\{v\in \mathbb{R}^n \mid (v,-1)\in \hat N_{{\rm epi}(\varphi)}(x,\varphi(x))\},\quad\forall x\in {\rm dom} (\varphi)$$
and
$$\partial\varphi(x):=\{v\in \mathbb{R}^n \mid (v,-1)\in N_{{\rm epi}(\varphi)}(x,\varphi(x))\},\quad\forall x\in {\rm dom} (\varphi).$$
Recall that the indicator function of $S$ is defined as
\begin{equation*}
\delta_S(x):=\left\{\begin{matrix}
   0,&{if\,x\in S},\\ +\infty,&{otherwise}.
\end{matrix}\right.
\end{equation*}
It is well known that
$$\hat N_S(x)=\hat \partial \delta_S(x)\,\,{\rm and}\,\, N_S(x)=\partial \delta_S(x),\quad\forall x\in S.$$

Now we are ready to review two kinds of second-order tangent sets which play a fundamental role in the second-order analysis later on.

\begin{definition}[Second-Order Tangent Sets \cite{shap,pen98}]\label{t1}
Let $S\subseteq \mathbb{R}^n$, $\bar x\in  S$ and $d \in T_S(\bar x)$ be given.

(i) The (outer) second-order tangent set to $S$ at $\bar x$ in direction $d$ is defined by
\begin{equation*}
T_S^2(\bar x;d):=\{w\in \mathbb{R}^n \mid \exists t_k\downarrow 0,w_k\rightarrow w\,\,{\rm such\,that}\,\,
\bar x +t_k d+\frac{1}{2}t_k^2w_k\in S\}.
\end{equation*}

(ii) The asymptotic second-order tangent cone to $S$ at $\bar x$ in direction $d$ is defined by
\begin{equation*}
\aligned
T_S^{''}(\bar x;d):=&\{w\in \mathbb{R}^n \mid \exists (t_k,r_k)\downarrow (0,0),w_k\rightarrow w\,\,{\rm such\,that}\\
&t_k/r_k\rightarrow 0,\bar x +t_k d+\frac{1}{2}t_kr_kw_k\in S \}.
\endaligned
\end{equation*}
\end{definition}

It is clear from the above definitions that both sets $T_S^2(\bar x;d)$ and $T_S^{''}(\bar x;d)$ are closed. The asymptotic second-order tangent cone was first introduced by Penot \cite{pen98} in the study of optimality conditions for scalar optimization, the original definition of which is slightly different while being equivalent to the one in Definition\eqref{t1}(ii). It is known that unlike the first-order tangent cone which is always a nonempty cone, the second-order tangent set may not be a cone and it may be empty (see, e.g., Bonnans and Shapiro \cite[Example 3.29]{shap}). However it can be easily verified that   $T_S^{''}(\bar x;d)$ is indeed a cone.

In general, both $T_S^{2}(\bar x;d)$ and $T_S^{''}(\bar x;d)$ are contained in 
{${\rm cl}([{\rm cone} [{\rm cone}(S-\bar x)-d]])$} and when $S$ is convex, $T_S^{''}(\bar x;d)={{\rm cl}([{\rm cone} [{\rm cone}(S-\bar x)-d]])}$ and $T_S^{2}(\bar x;d)\subseteq T_S^{''}(\bar x;d)$, which indicates that the asymptotic second-order tangent cone {may be} larger than the second-order tangent set.  For more details, please refer to \cite{xde} and the references therein.

In \cite[Proposition 2.1]{pen98}, it was  stated that for any $d\in T_S(\bar x)$, the sets $T_S^{2}(\bar x;d)$ and $T_S^{''}(\bar x;d)$   cannot be empty simultaneously. However, the set $T_S^{''}(\bar x;d)$  is a cone and, hence, it 
cannot be empty since $0\in T_S^{''}(\bar x;d)$ always.  Indeed, a slight modification of the proof of \cite[Proposition 2.1]{pen98} shows that the sets $T_S^{2}(\bar x;d)$ and $T_S^{''}(\bar x;d)\setminus \{0\}$   cannot be empty simultaneously. For completeness we revise the statements in \cite[Proposition 2.1]{pen98} and provide a proof here.
\begin{proposition}Let $S \subseteq \mathbb{R}^n$ be a closed set, $\bar x\in S$ and $d\in T_S(\bar x)$. Then,
\begin{equation*}\label{eqn22}
T_S^{2}(\bar x;d) \cup (T_S^{''}(\bar x;d)\backslash \{0\}) \neq \emptyset.\end{equation*}
\end{proposition}
\beginproof
If $d=0$, one has $T_S^{2}(\bar x;d)=T_S^{''}(\bar x;d)=T_S(\bar x)$ and the statement is obviously true. We now assume $d\not =0$ and 
 $T_S^2(\bar x;d) =\emptyset$.  Since $d \in T_S(\bar x)$,  there exist $t_k\downarrow 0,d_k\rightarrow d$ such that
$\bar x +t_k d_k\in S$.
Let $x_k:=\bar x +t_k d_k$, $\tilde{t}_k:=\frac{t_k\|d_k\|}{\|d\|}$ and $v_k:=\frac{x_k-\bar x-\tilde{t}_kd}{\frac{1}{2}\tilde{t}_k^2}$. Clearly, $\tilde{t}_k\downarrow 0$.
If $\{v_k\}$ is bounded, without loss of generality, we may assume that $v_k\rightarrow v$.
Since $x_k=\bar x +\tilde{t}_k d+\frac{1}{2}\tilde{t}_k^2v_k\in S$, we have $v\in T_S^2(\bar x;d)$, which is a contradiction to our assumption that $T_S^2(\bar x;d)=\emptyset$.
Therefore, $\{v_k\}$ is unbounded. 
Let $w_k:=\frac{v_k}{\|v_k\|}$. Without loss of generality, assume that $w_k\rightarrow w$, where $\|w\|=1$.
Let $r_k:=\tilde{t}_k\|v_k\|$. Then one has $\frac{\tilde{t}_k}{r_k}=\frac{1}{\|v_k\|}\rightarrow 0$, $r_k \downarrow 0$   
 (according to the definition of $v_k$ and $\frac{x_k-\bar x}{\tilde{t}_k}\rightarrow d$), 
and 
$$x_k=\bar x +\tilde{t}_k d+\frac{1}{2}\tilde{t}_k^2v_k=\bar x +\tilde{t}_k d+\frac{1}{2}\tilde{t}_kr_kw_k\in S.$$
This show that $w\in T_S^{''}(\bar x;d)$. Since $w\not =0$, the proof is completed.
\endproof

It is easy to obtain the following results from definition.
Given $S\subseteq \mathbb{R}^n, \bar x\in S$ and $d\in T_S(\bar x)$. Then we have
\begin{equation}\label{qctt}
T_S^2(\bar x;td)=t^2T_S^2(\bar x;d)\ \  \text{and}\ \  T_S^{''}(\bar x;td)={T_S^{''}(\bar x;d)}
\end{equation}
hold for all positive $t$.

We now derive some properties of first- and second-order tangent sets, which will be useful later on.

\begin{proposition}\label{ppzz}
Let $S \subseteq \mathbb{R}^n$ be a closed set, $\bar x\in S$ and $d\in T_S(\bar x)$. Then,
$$T_S^2(\bar x;d)+\hat T_S(\bar x;d)=T_S^2(\bar x;d)\ \  \text{and}\ \  T_S^{''}(\bar x;d)+\hat T_S(\bar x;d)=T_S^{''}(\bar x;d).$$
\end{proposition}
\beginproof
The first equality is clear from \cite[Proposition 1]{ye22}. It suffices to show that the second equality holds. The inclusion $\supseteq$ in second equation is trivial because $0\in \hat T_S(\bar x;d)$. For the converse inclusion $\subseteq$, consider $w\in T_S^{''}(\bar x;d)$ and $v\in \hat T_S(\bar x;d)$, and we prove that $w+v\in T_S^{''}(\bar x;d)$ by contradiction. Assume to the contrary that $w+v\not\in T_S^{''}(\bar x;d)$.
By definition, there exists $\epsilon>0$ and sequences $(t_k,r_k)\downarrow (0,0),w_k\rightarrow w$ such that
$t_k/r_k\rightarrow 0$ with $\bar x +t_k d+\frac{1}{2}t_kr_kw_k\in S$
and $$d\left(\bar x +t_k d+\frac{1}{2}t_kr_k(w_k+v),S\right) > \frac{1}{2}t_kr_k\epsilon$$
for {sufficiently large} $k\in \mathbb{N}$.
Denote by $x_k:=\bar x +t_k d+\frac{1}{2}t_kr_kw_k$ and $\tau_k:=\frac{1}{2}t_kr_k$, the preceding inequality is equivalent to saying that
$(x_k+\tau_kB_\epsilon(v))\cap S=\emptyset$. Then following the lines of the arguments in \cite[Proposition 1]{ye22}, we deduce that $v\not\in \hat T_S(\bar x;d)$, which is a contradiction.
\endproof

The following estimation on limiting normal cone to second-order tangent sets will be used in our development of directional second-order optimality conditions.
\begin{lemma}\label{xxzz}
Let $S$  be a closed subset in $\mathbb{R}^n$, $\bar x\in S$ and $d\in T_S(\bar x)$. Then,
$$N_{T_S^2(\bar x;d)}(w)\subseteq N_S(\bar x;d),\   \forall w\in T_S^2(\bar x;d), \mbox{ and } N_{T_S^{''}(\bar x;d)}(w)\subseteq N_S(\bar x;d),\  \forall w\in T_S^{''}(\bar x;d).$$
\end{lemma}
\beginproof
The first inclusion follows from \cite[Lemma 3]{ye22}. It suffices to show that the second inclusion holds. If $T_S^{''}(\bar x;d)=\emptyset$, the {assertion} is trivial. Otherwise, we first claim that
\begin{equation}\label{hhh}
\hat N_{T_S^{''}(\bar x;d)}(w)\subseteq N_S(\bar x;d),\quad \forall w\in T_S^{''}(\bar x;d).
\end{equation}
To this end, we pick any $w\in T_S^{''}(\bar x;d)$, which indicates the existence of sequences $(t_k,r_k)\downarrow (0,0),w_k\rightarrow w$ such that
$t_k/r_k\rightarrow 0$ with $\bar x +t_k d+\frac{1}{2}t_kr_kw_k\in S$ for all $k\in \mathbb{N}$.
Let {$T_k:=\{w'\in \mathbb{R}^n \mid \bar x +t_k d+\frac{1}{2}t_kr_kw'\in S\}$}. It is easy to observe that $T_k$ is closed, $w_k\in T_k$ and $w\in \limsup_{k\rightarrow \infty}T_k\subset T_S^{''}(\bar x;d)$. Then,
$$\hat N_{T_S^{''}(\bar x;d)}(w)\subseteq \hat N_{\limsup\limits_{k\rightarrow \infty}T_k}(w)\subseteq N_{\limsup\limits_{k\rightarrow \infty}T_k}(w).$$
For any $v\in \hat N_{T_S^{''}(\bar x;d)}(w)$, one has $v\in N_{\limsup\limits_{k\rightarrow \infty}T_k}(w)$.
It follows from Rockafellar and Wets \cite[Exercise 6.18]{rock} that there exist $w_k'\in T_k$ and $v_k\in \hat N_{T_k}(w_k')$, such that
$w_k'\rightarrow w$ and $v_k\rightarrow v$ (taking subsequences if necessary). Since $v_k\in \hat N_{T_k}(w_k')=\hat N_S(\bar x +t_k d+\frac{1}{2}t_kr_kw_k')$,
we have $v\in N_S(\bar x;d)$. Hence \eqref{hhh} holds.
Since $N_S(\bar x;d)$ is closed, it follows from \eqref{hhh} that $N_{T_S^{''}(\bar x;d)}(w)=\limsup\limits_{w'\rightarrow w}\hat N_{T_S^{''}(\bar x;d)}(w')\subseteq N_S(\bar x;d)$.
\endproof

{Note that both inclusions in Lemma \ref{xxzz} may be strict, see \cite[Example 1]{ye22} and the following example.} 
\begin{example}\label{ff1j}
Consider $S=\{(x_1,x_2)\in \mathbb{R}^2:x_1^2-x_2^3=0\}$ and $\bar x=(0,0)$. It is easy to calculate that $T_S(\bar x)=\mathbb{R}_+(0,1)$, where $\mathbb{R}_+(0,1)$ represents the cone generated by the vector $d=(0,1)$. Then we have $T_S^{''}(\bar x;d)=\mathbb{R}^2$ and $T_S^{2}(\bar x;d)=\emptyset$. 
Hence $N_{T_S^{''}(\bar x;d)}(w)=\{(0,0)\}$ for all $w\in \mathbb{R}^2$, while $N_S(\bar x;d)=\mathbb{R}\times \{0\}$, which implies that $N_{T_S^{''}(\bar x;d)}(w)$ is a proper subset of
$N_S(\bar x;d)$.
\end{example}

\subsection{Second-Order Tangent Set and Asymptotic Second-Order Tangent Cone to Feasible Set under Directional MSCQ}

In this subsection, we provide some auxiliary facts regarding the second-order tangent set and the asymptotic second-order tangent cone to the feasible set $C:=g^{-1}(K)$. The preliminary results are obtained under the validity of directional MSCQ of the constraint system $g(x)\in K$.

To introduce the property of directional MSCQ, we need to review the concept of directional neighborhood.
Let $d\in \mathbb{R}^n$ be a given direction and $\rho,\delta>0$, the directional neighborhood of direction $d$ is defined as follows:
\begin{equation*}
\aligned
V_{\rho,\delta}(d):&=\left\{w\in \delta B_{\mathbb{R}^n} \mid \left\| \|d\|w-\|w\|d\right\|\leq \rho\|w\|\|d\|\right\}\\
&=\left\{\begin{matrix}
   \delta B_{\mathbb{R}^n},&{if\,d=0},\\
   \left\{w\in \delta B_{\mathbb{R}^n}\backslash\{0\} \mid \left\|\frac{w}{\|w\|}-\frac{d}{\|d\|}\right\|\leq \rho\right\}\cup \{0\},&{if\,d\not=0}.
\end{matrix}\right.
\endaligned
\end{equation*}
Note that the directional neighborhood is not a ball except for the case of $d=0$. In general, when $d\in \mathbb{R}^n\backslash\{0\}$, the directional neighborhood is a section of the classical neighborhood.  Furthermore, for $\bar x\in \mathbb{R}^n$, if there are sequences $\{x_k\}\subseteq \mathbb{R}^n$ and $\{t_k\}\downarrow 0$ such that $\frac{x_k-\bar x}{t_k}\rightarrow d$, then
$
x_k\in \bar x+V_{\rho,\delta}(d),\,{\rm for\, sufficiently\, large}\,\, k.
$
Indeed, the aforementioned relation holds due to the fact that $\left\|\frac{x_k-\bar x}{t_k}\right\|\rightarrow \|d\|$
and
\begin{equation}\label{xkk}
\frac{x_k-\bar x}{\|x_k-\bar x\|}=\frac{x_k-\bar x}{t_k}\frac{t_k}{\|x_k-\bar x\|}\rightarrow\frac{d}{\|d\|}.
\end{equation}

Let $M:\mathbb{R}^n\rightrightarrows \mathbb{R}^m$ be a given multifunction and $(\bar x,\bar y)\in {\rm gph} M$. $M$ is said to be metrically subregular at $(\bar x,\bar y)$ in direction $d\in \mathbb{R}^n$ if there are positive numbers $\rho,\delta,\kappa>0$ such that
\begin{equation*}
d(x,M^{-1}(\bar y))\leq \kappa d(\bar y,M(x)),\quad \forall x\in \bar x+V_{\rho,\delta}(d).
\end{equation*}
The infimum of $\kappa$ over all the combinations $\rho,\delta$ and $\kappa$ satisfying the above relation is called the modulus of directional metric subregularity. In the case of $d=0$, we simply say that $M$ is metrically subregular at $(\bar x,\bar y)$.

For the constraint system $g(x) \in K$, we say that the metric subregularity constraint qualification (MSCQ) holds at $\bar x$ in direction $d$ {(with modulus $\kappa$)} if the set-valued mapping $M(x):=g(x)-K$ is metrically subregular at $(\bar x,0)$ in direction $d$ {(with modulus $\kappa$)}.

The following inclusion provides an upper estimate for the directional normals to the feasible set $C:=g^{-1}(K)$, which reduces to \cite[Theorem 3]{gfer17} when {direction $d=0$ is considered.}

\begin{lemma}\label{fxsb}
Let $\bar x\in C:=g^{-1}(K)$ and assume that MSCQ holds at $\bar x$ in direction $d\in \mathbb{R}^n$ for the constraint system $g(x)\in K$ with modulus $\kappa$. Then,
\begin{equation}\label{bbaa}
N_{g^{-1}(K)}(\bar x;d)\subseteq \{v\in \mathbb{R}^n \mid \exists \lambda\in N_K(g(\bar x);\nabla g(\bar x)d)\cap \kappa \|v\|B_{\mathbb{R}^m}\,\,{\rm with}\,\,v=\nabla g(\bar x)^T\lambda\}.
\end{equation}
\end{lemma}
\beginproof
Note that $g^{-1}(K)=M^{-1}(0)$ with $M(x):=g(x)-K$, by
applying \cite[Proposition 4.1]{gfer16}, we obtain that
$$N_{g^{-1}(K)}(\bar x;d)\subseteq \{v\in \mathbb{R}^n \mid \exists \lambda\in\kappa \|v\|B_{\mathbb{R}^m}\,\,{\rm with}\,\,(v,-\lambda)\in N_{{\rm gph}(M)}((\bar x,0);(d,0))\}.$$
Because ${\rm gph} M=\{(x,y) \mid g(x)-y\in K\}$, it follows from  \cite[Corollary 3.2]{ben} that
$$N_{{\rm gph}(M)}((\bar x,0);(d,0)) = \{(\nabla g(\bar x)^T\lambda,-\lambda) \mid \lambda\in N_K(g(\bar x);\nabla g(\bar x)d)\},$$
which yields \eqref{bbaa}.
\endproof

Now we are ready to present some basic facts about the second-order tangent set $T_C^2(\bar x;d)$ and the asymptotic second-order tangent cone $T_C^{''}(\bar x;d)$ to the feasible set $C:=g^{-1}(K)$ under directional MSCQ, which are useful in developing our main results later on.
\begin{proposition}\label{qc}
Let $\bar x\in C:=g^{-1}(K)$ and $d\in T_C(\bar x)$ be given. Then
\begin{eqnarray}
T_C^2(\bar x;d)&\subseteq & \{w\in \mathbb{R}^n \mid \nabla g(\bar x)w+\nabla^2 g(\bar x)(d,d)\in T_K^2(g(\bar x);\nabla g(\bar x)d)\}, \label{2cb}\\
T_C^{''}(\bar x;d)&\subseteq & \{w\in \mathbb{R}^n \mid \nabla g(\bar x)w\in T_K^{''}(g(\bar x);\nabla g(\bar x)d)\}.\label{2cb1}
\end{eqnarray}
If, in addition, MSCQ holds at $\bar x$ in direction $d\in \mathbb{R}^n$ for the constraint system $g(x)\in K$ with modulus $\kappa$, then inclusions \eqref{2cb} and \eqref{2cb1} hold as {equalities and the estimates}
\begin{eqnarray}\label{2cbb}
d(w,T_C^2(\bar x;d))&\leq& \kappa d(\nabla g(\bar x)w+\nabla^2 g(\bar x)(d,d), T_K^2(g(\bar x);\nabla g(\bar x)d)),
\\
\label{2cb1b}
d(w,T_C^{''}(\bar x;d))&\leq &\kappa d(\nabla g(\bar x)w, T_K^{''}(g(\bar x);\nabla g(\bar x)d))
\end{eqnarray}
hold for all $w\in \mathbb{R}^n$.
\end{proposition}
\beginproof
The results regarding the second-order tangent set $T_C^2(\bar x;d)$ follow from \cite[Proposition 5]{ye22}. Hence it suffices to prove the relations for the asymptotic second-order tangent cone, i.e., to show that \eqref{2cb1} and \eqref{2cb1b} hold true. To this end, we pick any $w\in T_C^{''}(\bar x,d)$, then there exist $(t_k,r_k)\downarrow (0,0)$ and $w_k\rightarrow w$ such that
$t_k/r_k\rightarrow 0$ and $x_k:=\bar x +t_k d+\frac{1}{2}t_kr_kw_k\in C$ for all $k\in \mathbb{N}$.
Then for all $ k\in \mathbb{N}$,
\begin{equation*}
\aligned
{K \ni }\  g(x_k)&=g(\bar x)+t_k\nabla g(\bar x)d+\frac{1}{2}t_kr_k\nabla g(\bar x)w_k+\frac{1}{2}t_k^2\nabla^2 g(\bar x)(d,d)+o(t_k^2)\\
&=g(\bar x)+t_k\nabla g(\bar x)d+\frac{1}{2}t_kr_k\left(\nabla g(\bar x)w_k+\frac{t_k}{r_k}\nabla^2 g(\bar x)(d,d)+\frac{o(t_k^2)}{\frac{1}{2}t_kr_k}\right).
\endaligned
\end{equation*}
Since $t_k/r_k\rightarrow 0$ and $w_k\rightarrow w$, one has
$\nabla g(\bar x)w_k+\frac{t_k}{r_k}\nabla^2 g(\bar x)(d,d)+\frac{o(t_k^2)}{\frac{1}{2}t_kr_k}\rightarrow \nabla g(\bar x)w$, which implies that
$\nabla g(\bar x)w\in T_K^{''}(g(\bar x);\nabla g(\bar x)d)$. Therefore, \eqref{2cb1} holds true.

Next, assume that the set-valued map $M(x):=g(x)-K$  is metrically subregular at $(\bar x,0)$ in direction $d$ with modulus $\kappa$. We aim to prove that estimate \eqref{2cb1b} holds which automatically {indicates that \eqref{2cb1} holds as equality.}
Pick any $\varepsilon>0$ and $\kappa'>\kappa$. Let $w\in \mathbb{R}^n$ be fixed, then there exists $u\in T_K^{''}(g(\bar x);\nabla g(\bar x)d)$
such that
\begin{equation}\label{bqq}
\|\nabla g(\bar x)w-u\|<d(\nabla g(\bar x)w, T_K^{''}(g(\bar x);\nabla g(\bar x)d))+\varepsilon.
\end{equation}
{Then} there exist sequences $(t_k,r_k) \downarrow (0,0)$ and $u_k \rightarrow u$ such that
$t_k/r_k\rightarrow 0$ and $g(\bar x) +t_k \nabla g(\bar x)d+\frac{1}{2}t_kr_ku_k\in K$ for all $k\in \mathbb{N}$.
Let $\tilde{x}_k:=\bar x +t_k d+\frac{1}{2}t_kr_kw$, then $\frac{\tilde{x}_k-\bar x}{t_k}=d+\frac{1}{2}r_kw\rightarrow d$.
So, for sufficiently large $k$, it follows from the directional metric subregularity of $M$ that
\begin{equation*}
\aligned
d(\tilde{x}_k,C)&\leq \kappa' d(g(\tilde{x}_k),K)\\
&\leq \kappa' \|g(\tilde{x}_k)-g(\bar x) -t_k \nabla g(\bar x)d-\frac{1}{2}t_kr_ku_k\|\\
&\leq \frac{\kappa' t_kr_k}{2}\left\|\nabla g(\bar x)w+\frac{t_k}{r_k}\nabla^2 g(\bar x)(d,d)-u_k\right\|+o(t_k^2).
\endaligned
\end{equation*}
Then
$$d\left(w,\frac{C-\bar x-t_k d}{\frac{1}{2}t_kr_k}\right)\leq \kappa'\left\|\nabla g(\bar x)w+\frac{t_k}{r_k}\nabla^2 g(\bar x)(d,d)-u_k\right\|+\frac{o(t_k^2)}{\frac{1}{2}t_kr_k},$$
which implies that there exists $x_k'\in C$ such that
\begin{equation}\label{xkp}
\|w-w_k\|\leq \kappa'\left\|\nabla g(\bar x)w+\frac{t_k}{r_k}\nabla^2 g(\bar x)(d,d)-u_k\right\|+4\frac{o(t_k^2)}{t_kr_k},
\end{equation}
where $w_k:= \frac{x_k'-\bar x-t_k d}{\frac{1}{2}t_kr_k}$. Since the {right-hand} side of \eqref{xkp} converges to $\kappa'\left\|\nabla g(\bar x)w-u\right\|$, we know that $\{w_k\}$ is bounded. Without loss of generality, we assume that $w_k\rightarrow w'$, then $w'\in T_C^{''}(\bar x;d)$.
Passing to the limit in \eqref{xkp}, we arrive at
$$d(w,T_C^{''}(\bar x;d))\leq\|w-w'\|\leq \kappa'\left\|\nabla g(\bar x)w-u\right\|\leq\kappa'd(\nabla g(\bar x)w, T_K^{''}(g(\bar x);\nabla g(\bar x)d))+ 
{\kappa'}\varepsilon.$$
Because $\kappa'$ can be chosen arbitrarily close to $\kappa$ and $\varepsilon$ can be chosen  arbitrarily small, we know that \eqref{2cb1b} holds true.
 From \eqref{2cb1b}, we may conclude that
$$T_C^{''}(\bar x;d)\supseteq \{w\in \mathbb{R}^n:\nabla g(\bar x)w\in T_K^{''}(g(\bar x);\nabla g(\bar x)d)\}, $$
which indicates that inclusion \eqref{2cb1} holds as equality.
\endproof

In the following proposition, we provide upper estimations for limiting normals to second order tangent set and asymptotic second-order tangent cone under directional MSCQ, respectively.
\begin{proposition}\label{sub}
Let $\bar x\in C:=g^{-1}(K)$ and {$d\in T_C(\bar x)$} be given. Assume that MSCQ holds at $\bar x$ in direction $d\in \mathbb{R}^n$ for the constraint system $g(x)\in K$ with modulus $\kappa$. Then,
\begin{equation*}
\aligned
N_{T_C^{''}(\bar x;d)}(v)\subseteq \{z&\in \mathbb{R}^n \mid \exists \lambda\in N_{T_K^{''}(g(\bar x);\nabla g(\bar x)d)}(\nabla g(\bar x)v)\cap \kappa \|z\|B_{\mathbb{R}^m}\\&{\rm with}\,\,z=\nabla g(\bar x)^T\lambda\},\quad \forall v\in T_C^{''}(\bar x,d)
\endaligned
\end{equation*}
and
\begin{equation*}
\aligned
N_{T_C^{2}(\bar x;d)}(w)\subseteq \{z&\in \mathbb{R}^n \mid \exists \lambda\in N_{T_K^{2}(g(\bar x);\nabla g(\bar x)d)}(\nabla g(\bar x)w+\nabla^2 g(\bar x)(d,d))\cap \kappa \|z\|B_{\mathbb{R}^m}\\&{\rm with}\,\,z=\nabla g(\bar x)^T\lambda\},\quad \forall w\in T_C^{2}(\bar x;d).
\endaligned
\end{equation*}
\end{proposition}

\beginproof
According to the assumption, for any $v\in T_C^{''}(\bar x,d)$, it follows from Proposition \ref{qc} that the mapping $w\mapsto \nabla g(\bar x)w- T_K^{''}(g(\bar x);\nabla g(\bar x)d)$ is metrically subregular at $(v,0)$ with modulus $\kappa$.
Then the first inclusion follows directly from \cite[Theorem 3]{gfer17} (it can also be obtained from \cite[Proposition 4.1]{gfer16} by choosing the zero direction). The second inclusion can be obtained similarly.
\endproof


\section{Optimality Conditions for Nonconvex Optimization Problem with an Abstract Set Constraint}
\label{sec:alg}
In this section we consider the following nonconvex optimization problem with an abstract set constraint:
\begin{equation}\label{opc}
\min f(x)\quad \text{s.t.}\quad x\in C,
\end{equation}
where $f:\mathbb{R}^n\rightarrow \mathbb{R}$ is twice continuously differentiable and $C$ is a closed subset of $\mathbb{R}^n$.
\begin{definition}[Directional Local Optimality] A point $\bar x\in C$ is said to be a local optimal solution of \eqref{opc} in direction $d\in \mathbb{R}^n$, if there exist positive numbers $\rho,\delta>0$ such that
\begin{equation}\label{los}
f(x)\geq f(\bar x),\,\,\,\forall x\in C\cap (\bar x+V_{\rho,\delta}(d)).
\end{equation}
\end{definition}

Note that the directional local optimal solution is only well-defined for directions $d\in T_C(\bar x)$. In fact, we claim that for an arbitary given $d\not\in T_C(\bar x)$, there exist $\rho,\delta>0$ such that $C\cap (\bar x+V_{\rho,\delta}(d))=\{\bar x\}$. If this is not the case, then
for each $k\in \mathbb{N}$ there exists $x_k\in C\cap (\bar x+V_{1/k,1/k}(d))\backslash\{\bar x\}$, that is, $x_k\in C, \|x_k-\bar x\|\leq\frac{1}{k}\rightarrow 0$ with $\left\|\frac{x_k-\bar x}{\|x_k-\bar x\|}-\frac{d}{\|d\|}\right\|\leq\frac{1}{k}\rightarrow 0$. We arrive at $\frac{d}{\|d\|}\in T_C(\bar x)$, which is a contradiction to our choice of $d\not\in T_C(\bar x)$.
Hence, $\bar x$ is always a local optimal solution of \eqref{opc} in each direction outside of $T_C(\bar x)$.
It worth to note that, if $\bar x$ is a local optimal solution of \eqref{opc} in all directions  $d\in T_C(\bar x)$, then $\bar x$ is a local optimal solution of \eqref{opc}.

For optimization problem with an abstract set constraint, there are well-developed theoretical results which reveals essential connections between variational geometry and optimality conditions. First-order necessary conditions for local optimality are provided in \cite[Theorem 6.12]{rock} in both primal and dual forms. In the following proposition, we improved the aforementioned results to the directional case which sharpened the necessary condition while weakened the optimality assumption. The obtained results are useful in developing our main results in Section \ref{3.2} and have its own independent interest.

\begin{proposition}[First-Order Necessary Optimality Condition]\label{op1}
Let $f:\mathbb{R}^n\rightarrow \mathbb{R}$ be continuously differentiable and $C\subset \mathbb{R}^n$ be closed. Suppose that $\bar x\in C, d\in T_C(\bar x)$ and $\bar x$ is a local optimal solution of \eqref{opc} in direction $d$.
Then,
\begin{enumerate}[{\rm (i)}]
\item $\nabla f(\bar x)d\geq 0$;
\item if $\nabla f(\bar x)d=0$, then $0\in\nabla f(\bar x)+N_C(\bar x;d)$.
\end{enumerate}
\end{proposition}
\beginproof
According to the assumptions, there exist $\rho,\delta>0$ such that \eqref{los} holds.
Since $d\in T_C(\bar x)$, there exist $t_k\downarrow 0$ and $d_k\rightarrow d$ such that
$x_k:=\bar x +t_k d_k\in C$ for all $k\in \mathbb{N}$.
Notice that $x_k\in C \cap (\bar x + V_{\rho,\delta}(d))$ for sufficiently large $k$ thanks to \eqref{xkk}, it then follows from \eqref{los} that
\begin{equation}\label{fxx}
f(x_k)=f(\bar x)+t_k\nabla f(\bar x)d_k+o(t_k)\geq f(\bar x).
\end{equation}
Dividing both sides of the above inequality by $t_k$ and passing to the limit, we conclude that (i) holds.

It remains to show the validity of (ii). Since $\frac{d_k}{\|d_k\|}-\frac{d}{\|d\|}\rightarrow 0$ and $x_k-\bar x=t_kd_k\rightarrow 0$, without loss of generality, we may assume that
$\left\|\frac{d_k}{\|d_k\|}-\frac{d}{\|d\|}\right\|<\frac{\rho}{2}$ and $\|x_k-\bar x\|<\frac{\delta}{2}$, for all $k\in \mathbb{N}$. For each $k$, let $r_k:=\frac{1}{4}\rho t_k\|d_k\|$, then, for any $x\in B_{r_k}(x_k)$, we have
\begin{equation*}
\aligned
\left\|\frac{x-\bar x}{\|x-\bar x\|}-\frac{d}{\|d\|}\right\|&\leq
\left\|\frac{x-\bar x}{\|x-\bar x\|}-\frac{x_k-\bar x}{\|x_k-\bar x\|}\right\|+\left\|\frac{d_k}{\|d_k\|}-\frac{d}{\|d\|}\right\|\\
&\leq \frac{2\|x-x_k\|}{\|x_k-\bar x\|}+\frac{\rho}{2}\leq\rho,
\endaligned
\end{equation*}
where the second inequality holds due to the fact that $\left\|\frac{u}{\|u\|}-\frac{v}{\|v\|}\right\|
\leq \frac{2\|u-v\|}{\|u\|}$ for all $u,v\in \mathbb{R}^n\backslash\{0\}$.
This shows that $B_{r_k}(x_k)\subset \bar x+V_{\rho,\delta}(d)$.
Again, by \eqref{los}, we have
$$f(x)\geq f(\bar x)\quad\forall x\in C\cap B_{r_k}(x_k).$$
Since $f$ is continuously differentiable, we can find $\varepsilon_k\in (0,1)$ such that $\varepsilon_k\rightarrow 0$ and
\begin{equation}\label{fai}
f(\bar x)\leq f(x_k)+\nabla f(x_k)(x-x_k)+\varepsilon_k\|x-x_k\|\quad\forall x\in C\cap B_{r_k}(x_k).
\end{equation}
Note that $\nabla f(\bar x)d=0$, by \eqref{fxx}, we have
$f(x_k)-f(\bar x)=t_k\nabla f(\bar x)d_k+o(t_k)=o(t_k).$
Without loss of generality, we assume that
\begin{equation}\label{faixx}
f(x_k)-f(\bar x)\leq \varepsilon_kt_k\quad \forall k\in \mathbb{N}.
\end{equation}
Define $\varphi_k:\mathbb{R}^n\rightarrow\mathbb{R}\cup\{+\infty\}$ as follows:
$$\varphi_k(x):=f(x_k)-f(\bar x)+\nabla f(x_k)(x-x_k)+\varepsilon_k\|x-x_k\|+\delta_{C\cap B_{r_k}(x_k)}(x)\quad\forall x\in\mathbb{R}^n.$$
It follows from \eqref{fai} that $\varphi_k(x_k)\leq \inf_{x\in\mathbb{R}^n}\varphi_k(x)+f(x_k)-f(\bar x)$.
By Ekeland's variational principle, there exists $u_k\in \mathbb{R}^n$ such that
$\|u_k-x_k\|\leq \sqrt{\varepsilon_k}r_k/4$, $\varphi_k(u_k)\leq \varphi_k(x_k)<+\infty$ and
\begin{equation}\label{ekf}
\varphi_k(u_k)\leq \varphi_k(x)+\frac{4}{\sqrt{\varepsilon_k}r_k}(f(x_k)-f(\bar x))\|x-u_k\|\quad\forall x \in \mathbb{R}^n.
\end{equation}
Then $u_k\in C\cap B_{r_k}(x_k)$.
By \eqref{ekf}, according to the first-order optimality conditions, we have
$0\in \hat{\partial}(\varphi_k+\frac{4}{\sqrt{\varepsilon_k}r_k}(f(x_k)-f(\bar x))\|\cdot-u_k\|)(u_k)$.
Then, by fuzzy sum rule of the subdifferential, there exists $u_k'\in B_{\sqrt{\varepsilon_k}r_k/4}(u_k)$ such that
\begin{equation}\label{cwf}
0\in \nabla f(x_k)+\left(\varepsilon_k+\frac{4}{\sqrt{\varepsilon_k}r_k}(f(x_k)-f(\bar x))+\frac{\sqrt{\varepsilon_k}r_k}{4}\right)B_{\mathbb{R}^n}+\hat N_{C\cap B_{r_k}(x_k)}(u_k').
\end{equation}
It is easy to see that
$\|u_k'-x_k\|\leq \|u_k'-u_k\|+\|u_k-x_k\|\leq \sqrt{\varepsilon_k}r_k/2<r_k$
and
\begin{equation*}
\aligned
\left\|\frac{u_k'-\bar x}{\|u_k'-\bar x\|}-\frac{d}{\|d\|}\right\|&\leq
\left\|\frac{u_k'-\bar x}{\|u_k'-\bar x\|}-\frac{x_k-\bar x}{\|x_k-\bar x\|}\right\|+\left\|\frac{d_k}{\|d_k\|}-\frac{d}{\|d\|}\right\|\\
&\leq \frac{2\|u_k'-x_k\|}{\|x_k-\bar x\|}+\left\|\frac{d_k}{\|d_k\|}-\frac{d}{\|d\|}\right\|\rightarrow 0.
\endaligned
\end{equation*}
This shows that $u_k'\in {\rm int}(B_{r_k}(x_k))$ and $d_k':=\frac{u_k'-\bar x}{t_k'}\rightarrow d$, where $t_k':=\|d\|^{-1}\|u_k'-\bar x\|$.
Then $\hat N_{C\cap B_{r_k}(x_k)}(u_k')=\hat N_{C\cap B_{r_k}(x_k)}(\bar x+t_k'd_k')=\hat N_{C}(\bar x+t_k'd_k')$.
It follows from \eqref{cwf} that, there exists $b_k\in B_{\mathbb{R}^n}$ such that
$$-\nabla f(x_k)+\left(\varepsilon_k+\frac{4}{\sqrt{\varepsilon_k}r_k}(f(x_k)-f(\bar x))+\frac{\sqrt{\varepsilon_k}r_k}{4}\right)b_k\in
\hat N_{C}(\bar x+t_k'd_k').$$
By \eqref{faixx}, we have
$$0\leq\frac{4}{\sqrt{\varepsilon_k}r_k}(f(x_k)-f(\bar x))\leq \frac{16\varepsilon_kt_k}{\sqrt{\varepsilon_k}\rho t_k\|d_k\|}
\rightarrow 0.$$
Since $\|b_k\|\leq 1,\nabla f(x_k)\rightarrow \nabla f(\bar x)$ and $\varepsilon_k\rightarrow 0$,
we conclude that (ii) holds.
\endproof

When the direction $d\neq 0$, the directional normal cone is smaller than the limiting normal cone. Hence Proposition~\ref{op1} obtains a stronger necessary optimality condition under weaker assumptions than the classical result \cite[Theorem 6.12]{rock}. Note that along with the condition $d\in T_C(\bar x)$ and Proposition \ref{op1} (i),  the condition of $\nabla f(\bar x)d=0$ in Proposition \ref{op1} (ii) means that $d$ must be a critical direction. The following example illustrates that, the condition of $\nabla f(\bar x)d=0$ in Proposition \ref{op1} (ii) cannot be dropped, even in the case when $C$ is polyhedral. It is worth to note that, when $C$ is convex, it follows from \eqref{ccc} that $N_C(\bar x;d)=N_C(\bar x)\cap \{d\}^\bot$. Therefore, if $-\nabla f(\bar x)\in N_C(\bar x;d)$, we must have $\nabla f(\bar x)d=0$.
\begin{example}\label{1j}
Let $f:\mathbb{R}^2\rightarrow\mathbb{R}$ be defined as $f(x_1,x_2)=(x_1+1)^2$ for all $(x_1,x_2)\in \mathbb{R}^2$. Consider $C=\mathbb{R}^2_+, \bar x=(0,0)$ and $d=(1,1)$. Then $\bar x$ is a local optimal solution of \eqref{opc} in direction $d$. But, in this case,
$\nabla f(\bar x)d\not=0$ and $-\nabla f(\bar x)\not\in N_C(\bar x;d)=\{0\}$.
\end{example}

Since the asymptotic second-order tangent cone $T_C^{''}(\bar x;d)$ and the outer second-order tangent set $T_C^{2}(\bar x;d)$ cannot be empty simultateously, it is more appropriate to consider second-order optimality conditions that involves both sets. The following proposition is an enhancement to \cite[Corollary 1.3]{pen94} and \cite[Corollary 3.2]{pen98} in which the author provided second-order necessary conditions for problem \eqref{opc} to obtain a local minimum with no direction specified. It establishes exact correspondence between the direction along which problem \eqref{opc} obtains a local solution and the direction involved in the presentation of the optimality conditions.
\begin{proposition}[Primary Second-Order Necessary Optimality Condition]\label{op2}
Let $d\in T_C(\bar x)$ and $\bar x$ be a local optimal solution of \eqref{opc} in direction $d\in \mathbb{R}^n$ with $\nabla f(\bar x)d=0$.
Then,
\begin{enumerate}[{\rm (i)}]
\item $\nabla f(\bar x)v \geq 0,\quad\forall v\in T_C^{''}(\bar x,d)$;
\item $\nabla f(\bar x)(w)+\nabla^2 f(\bar x)(d,d)\geq0,\quad\forall w\in T_C^2(\bar x;d)$.
\end{enumerate}
\end{proposition}

\beginproof
Since $\bar x$ is a local optimal solution of \eqref{opc} in direction $d\in \mathbb{R}^n$, there exist $\rho,\delta>0$ such that \eqref{los} holds.
Pick any $v\in T_C^{''}(\bar x;d)$, then there exist $(t_k,r_k)\downarrow (0,0)$ and $v_k\rightarrow v$ such that
$t_k/r_k\rightarrow 0$ and $x_k:=\bar x +t_k d+\frac{1}{2}t_kr_kv_k\in C$ for all $k\in \mathbb{N}$.
Since $\frac{x_k-\bar x}{t_k}=d+r_kv_k\rightarrow d$, then we have \eqref{xkk} holds.
It follows from \eqref{los} that, for sufficiently large $k$,
\begin{equation*}
f(x_k)=f(\bar x)+t_k\nabla f(\bar x)d+\frac{1}{2}t_kr_k\nabla f(\bar x)v_k+\frac{1}{2}t_k^2\nabla^2 f(\bar x)\left(d,d\right)+o(t_k^2)\geq f(\bar x).
\end{equation*}
Note that $\nabla f(\bar x)d=0$, one has
\begin{equation*}
\nabla f(\bar x)v_k+\frac{t_k}{r_k}\nabla^2 f(\bar x)\left(d,d\right)+\frac{o(t_k^2)}{\frac{1}{2}t_kr_k}\geq 0.
\end{equation*}
Since $t_k/r_k\rightarrow 0$, passing to the limit in the above inequality, we arrive at (i).

Pick any $w\in T_C^2(\bar x;d)$. Then there exist $t_k\downarrow 0$ and $w_k\rightarrow w$ such that
$\tilde{x}_k:=\bar x +t_k d+\frac{1}{2}t_k^2w_k\in C$ for all $k\in \mathbb{N}$.
Similarly,  we know that for sufficiently large $k$, $\tilde{x}_k\in C\cap (\bar x+V_{\rho,\delta}(d))$ and
\begin{equation*}
f(x_k)=f(\bar x)+t_k\nabla f(\bar x)d+\frac{1}{2}t_k^2(\nabla f(\bar x)w_k+\nabla^2 f(\bar x)(d,d))+o(t_k^2)\geq f(\bar x).
\end{equation*}
Dividing both sides of the above inequality by $\frac{1}{2}t_k^2$ and passing to the limit, we have (ii) holds.
\endproof

Assertion (ii) in Proposition \ref{op2} shows that $\inf_{w\in T_C^2(\bar x;d)}(\nabla f(\bar x)(w)+\nabla^2 f(\bar x)(d,d))$ $\geq 0$ when $\bar x$ is a local optimal solution of \eqref{opc} in direction $d$. In the following example, we illustrate that this infimum can be strictly positive. 

\begin{example}\label{fj2}
Let $f:\mathbb{R}^2\rightarrow\mathbb{R}$ be defined as $f(x_1,x_2)=x_2^2$ for all $(x_1,x_2)\in \mathbb{R}^2$. Consider $C=\mathbb{R}^2_+, \bar x=(0,0)$ and $d=(0,1)$, then $\bar x$ is a local optimal solution of \eqref{opc} in direction $d$. It is easy to see that
$\nabla f(\bar x)=(0,0),\nabla^2 f(\bar x)=\left(
\begin{array}{cc}
 0 & 0\\
 0 & 2\\
\end{array}
\right)
$ and $T_C^2(\bar x;d)=\mathbb{R}_+\times \mathbb{R}$. Hence, we have $\inf_{w\in T_C^2(\bar x;d)}(\nabla f(\bar x)(w)+\nabla^2 f(\bar x)(d,d))=2> 0$.
\end{example}


\section{Optimality Conditions for  Problem (\ref{op})}
\label{3.2}

In this section, we aim at developing new necessary and sufficient optimality conditions for problem~\eqref{op} to characterize a local optimal solution in direction $d$. 
Throughout this section we {denote by  $C:=g^{-1}(K)$, 
 $\Theta:= \nabla g(\bar x)(T_C^{''}(\bar x; d))$ and $\Omega:=\nabla g(\bar x)(T_C^2(\bar x;d))+\nabla^2 g(\bar x)(d,d)$}.}

Since  problem~\eqref{op} is problem~\eqref{opc} with $C:=g^{-1}(K)$, in the following proposition, by Proposition  \ref{op1}, we provide the first-order necessary optimality condition for problem~\eqref{op} at a local directional minimizer under the directional MSCQ for the constraint system $g(x)\in K$. Note that one always has
$$ d\in T_C(\bar x) \Longrightarrow  \nabla g(\bar x)d \in T_K(g(\bar{x})) $$ and the equivalence $\Leftrightarrow $ holds if the  directional MSCQ holds in direction $d$.

\begin{proposition}[First-Order Necessary Optimality Condition]\label{op22}
Let $f:\mathbb{R}^n\rightarrow \mathbb{R}$ and $g:\mathbb{R}^n\rightarrow \mathbb{R}^m$ be continuously differentiable, $d\in T_C(\bar x)$ and $\bar x$ be a local optimal solution of problem \eqref{op} in direction $d$. Suppose that $\nabla f(\bar x)d=0$ and MSCQ holds at $\bar x$ in direction $d$ for the constraint system $g(x)\in K$ with modulus $\kappa$. Then
 there exists $\lambda\in N_K(g(\bar x);\nabla g(\bar x)d)$ such that
\begin{equation}\label{oop}
\|\lambda\|\leq \kappa\|\nabla f(\bar x)\|\,\,{\rm and}\,\,\nabla_x L(\bar x,\lambda)=0.
\end{equation}
\end{proposition}

\beginproof
By Proposition \ref{op1}, we have $-\nabla f(\bar x)\in N_C(\bar x;d)$. At the hand of directional MSCQ, it follows from Lemma \ref{fxsb} that
\begin{eqnarray*}
-\nabla f(\bar x)\in N_C(\bar x;d)\subseteq \{v\in \mathbb{R}^n \mid &\exists \lambda\in N_K(g(\bar x);\nabla g(\bar x)d)\cap \kappa \|v\|B_{\mathbb{R}^m} \\&{\rm with}\,\,v=\nabla g(\bar x)^T\lambda\}.
\end{eqnarray*}
This implies that there exists $\lambda\in N_K(g(\bar x);\nabla g(\bar x)d)$ such that $-\nabla f(\bar x)=\nabla g(\bar x)^T\lambda$, i.e., $\nabla_x L(\bar x,\lambda)=0$, and
$\|\lambda\|\leq \kappa\|\nabla g(\bar x)^T\lambda\|=\kappa\|\nabla f(\bar x)\|$.
\endproof
{For convenience, for each feasible solution $\bar x$ and each $d\in T_C(\bar x)$, define the   M-multiplier set in direction $d$ as
$$\Lambda(\bar x;d):=\{\lambda\in \mathbb{R}^n \mid \lambda\in N_K(g(\bar x);\nabla g(\bar x)d)\,\,{\rm and}\,\,\nabla_x L(\bar x,\lambda)=0\}.$$
By Proposition \ref{op22}, the directional MSCQ ensures that $\Lambda(\bar x;d)$ for a local optimal solution $\bar x$ and a direction $d\in T_C(\bar x)\cap \{\nabla f(\bar x) \}^\perp$ is nonempty.}

We now prove a new second-order necessary optimality condition.
\begin{theorem}[Second-Order Necessary  Optimality Condition]\label{nss}
Let $d\in T_C(\bar x)$ and $\bar x$ be a local optimal solution of problem \eqref{op} in direction $d$ with $\nabla f(\bar x)d=0$.
Then, for any $\lambda\in \mathbb{R}^m$ with $\nabla_x L(\bar x,\lambda)=0$, one has
\begin{enumerate}[{\rm (i)}]
\item  $\sigma_\Theta(\lambda) \leq 0$;
\item $\nabla_{xx}^2L(\bar x,\lambda)(d,d)-\sigma_\Omega
(\lambda)= \alpha\geq0$, where
$$\alpha:=\inf_{w\in T_C^2(\bar x;d)}(\nabla f(\bar x)(w)+\nabla^2 f(\bar x)(d,d)).$$
\end{enumerate}
\end{theorem}

\beginproof
Since $\bar x$ is a local optimal solution of \eqref{op} in direction $d\in \mathbb{R}^n$, it follows from Proposition \ref{op2} that
$\alpha\geq 0$ and
\begin{equation}\label{i}
\nabla f(\bar x)u\geq 0,\quad\forall u\in T_C^{''}(\bar x;d).
\end{equation}
Pick any $\lambda$ such that $\nabla_x L(\bar x,\lambda)=0$. For any $v\in \nabla g(\bar x)(T_C^{''}(\bar x;d))$, there exists $u\in T_C^{''}(\bar x;d)$ such that $v=\nabla g(\bar x)u$.
It follows from \eqref{i} that
$$\langle\lambda,v\rangle=\langle\nabla g(\bar x)^T\lambda,u\rangle=-\nabla f(\bar x)u\leq 0,$$
which shows that (i) holds. We have
\begin{equation*}
\aligned
\sigma_{\Omega}(\lambda)&=\sup_{u\in \Omega}\langle\lambda,u\rangle\\
&=\sup_{w\in T_C^2(\bar x;d)}\langle\lambda,\nabla g(\bar x)w+\nabla^2 g(\bar x)(d,d)\rangle
\\&=\sup_{w\in T_C^2(\bar x;d)}\langle\nabla g(\bar x)^T\lambda,w\rangle+\langle\lambda,\nabla^2 g(\bar x)(d,d)\rangle
\\&=\sup_{w\in T_C^2(\bar x;d)}-(\nabla f(\bar x)w+\nabla^2 f(\bar x)(d,d))+\nabla^2 f(\bar x)(d,d)+\langle\lambda,\nabla^2 g(\bar x)(d,d)\rangle
\\&=-\inf_{w\in T_C^2(\bar x;d)}(\nabla f(\bar x)w+\nabla^2 f(\bar x)(d,d))+\nabla^2_{xx} L(\bar x,\lambda)(d,d)
\\&=-\alpha+\nabla^2_{xx} L(\bar x,\lambda)(d,d),
\endaligned
\end{equation*}
which implies that (ii) holds.
\endproof

In the following example, we show that even when MSCQ holds with $K$ being convex, the inclusion \eqref{strin} can be strict, and hence $\sigma_{\Omega}$ in Theorem~\ref{nss}(ii) cannot be replaced by {$\sigma_{T^2_K (g(\bar x);\nabla g(\bar x)d)}(\lambda)$.}

\begin{example}\label{fg3a}
Consider the simple constrained optimization problem
\begin{equation*}
\min f(x):=\frac{x^2}{2}\,\,{\rm s.t.}\,\,g(x):=(x^2,x)\in K:=\{(u,v)\in \mathbb{R}^2 \mid (u-1)^2+v^2\leq 1\}.
\end{equation*}
Let $\bar x=0\in \mathbb{R}$ and $d=1\in \mathbb{R}$.  Clearly, the feasible set computes as {$C:=g^{-1}(K)=[-1,1]$}, and thus $\bar x$ is a local optimal solution in direction $d$ and {it is straightforward to verify from the definition that the mapping $g(\cdot)-K$ is metrically subregular at $(\bar x,0)$ with $\kappa=1$ and we have $d(x,g^{-1}(K))\leq \kappa d(g(x),K)$ holds with equality.} By direct calculations we have $d\in T_C(\bar x)=\mathbb{R}_+$, $T_C^{''}(\bar x;d)=T_C^2(\bar x;d)=\mathbb{R}$, $\nabla g(\bar x)d=(0,1)$, $\nabla^2g(\bar x)(d,d)=(2,0)$
and $\Omega:=\nabla g(\bar x)(T_C^2(\bar x;d))+\nabla^2 g(\bar x)(d,d)=\{2\}\times \mathbb{R}$.
Utilizing {\cite[Propositions 13.13]{rock}}, it is not
difficult to show that
 $T^2_K (g(\bar x);\nabla g(\bar x)d)=\{(w_1,w_2)\in \mathbb{R}^2:w_1\geq 1\}$.
So, $\Omega\subsetneqq T^2_K (g(\bar x);\nabla g(\bar x)d)$.
Since $\nabla f(\bar x)=0$, Theorem \ref{nss} (i) holds. Pick any $\lambda=(\lambda_1,\lambda_2)$ satisfying $\nabla_x L(\bar x,\lambda)=0$, then $\nabla g(\bar x)^T\lambda=-\nabla f(\bar x)=0$, and then $\lambda_2=0$.
Direct calculation yields $\nabla^2_{xx}L(\bar x,\lambda)(d,d)-\sigma_\Omega(\lambda)=1+2\lambda_1-2\lambda_1>0$, which implies that Theorem \ref{nss} (ii) holds as well. However, for $\bar\lambda=(-2,0)$ which satisfies $\nabla_x L(\bar x,\bar\lambda)=0$, we have $\sigma_{T^2_K (g(\bar x);\nabla g(\bar x)d)}(\bar\lambda)=-2>-4=\sigma_\Omega(\bar\lambda)$ and
$$\nabla^2_{xx}L(\bar x,\bar\lambda)(d,d)-\sigma_{T^2_K (g(\bar x);\nabla g(\bar x)d)}(\bar\lambda)=-1<0<\nabla^2_{xx}L(\bar x,\bar\lambda)(d,d)-\sigma_\Omega(\bar\lambda).$$
\end{example}

In Theorem \ref{nss}, we only know that $\lambda$ satisfies the condition $\nabla_xL(\bar x,\lambda)=0$.  In the following result, we show that under the direcional MSCQ, the value of  $\lambda$ in  Theorem \ref{nss} is a directional M-multiplier with some additional properties.

\begin{theorem}\label{lmby}
Let $d\in T_C(\bar x)$ and $\bar x$ be a local optimal solution of problem \eqref{op} in direction $d$. Suppose that $\nabla f(\bar x)d=0$ and MSCQ holds at $\bar x$ in direction $d$ for the constraint system $g(x)\in K$. 
Then there exists $\lambda\in \Lambda(\bar x;d)$ such that $\sigma_\Theta(\lambda) \leq 0$. Moreover, $\lambda\in N_{T_K^{''}(g(\bar x);\nabla g(\bar x)d)}(0)$.
If $T^2_K (g(\bar x);\nabla g(\bar x)d)\not=\emptyset$, then {there exist $\lambda\in \Lambda(\bar x;d)$ such that  $$\nabla_{xx}^2L(\bar x,\lambda)(d,d)-\sigma_\Omega
(\lambda)\geq0. $$ Moreover there exists sequence $w_k\in T_C^2(\bar x;d)$ and $\lambda_k\in N_{T_K^{2}(g(\bar x);\nabla g(\bar x)d)}(\nabla g(\bar x)w_k+\nabla^2 g(\bar x)(d,d))$ such that $\nabla f(\bar x)w_k+\nabla^2 f(\bar x)(d,d)$ converges to $\alpha$ as defined in Theorem \ref{nss}(ii), $\langle\nabla g(\bar x)^T\lambda_k+\nabla f(\bar x),w_k\rangle\rightarrow 0$ and $\lambda_k$ converges to  $\lambda$.}

\end{theorem}

\beginproof
By Proposition \ref{op2} (i), we know that $-\nabla f(\bar x)\in (T_C^{''}(\bar x,d))^\circ$. Since $0\in T_C^{''}(\bar x;d)$, we have $(T_C^{''}(\bar x,d))^\circ\subseteq N_{T_C^{''}(\bar x,d)}(0)$, which implies that 
$$-\nabla f(\bar x)\in N_{T_C^{''}(\bar x,d)}(0).$$
{Since  MSCQ holds at $\bar x$ in direction $d$, by Proposition \ref{sub}, there exists $\lambda\in N_{T_K^{''}(g(\bar x);\nabla g(\bar x)d)}(0)$ such that $-\nabla f(\bar x)=\nabla g(\bar x)^T\lambda$.}
By Lemma \ref{xxzz}, we have $\lambda\in N_{T_K^{''}(g(\bar x);\nabla g(\bar x)d)}(0)\subseteq N_K(g(\bar x);\nabla g(\bar x)d)$, and then $\lambda\in \Lambda(\bar x;d)$. 
 It follows from Theorem \ref{nss}(i) that $\sigma_\Theta(\lambda) \leq 0$. 

Now assume that $T^2_K (g(\bar x);\nabla g(\bar x)d)\not=\emptyset$.  By Proposition \ref{op2} (ii), we have $\alpha=\inf_{w\in T_C^2(\bar x;d)}(\nabla f(\bar x)w+\nabla^2 f(\bar x)(d,d))\geq 0$. {From Proposition \ref{qc}, we know that $T_C^2(\bar x;d)$ is not empty.}
Let $\varphi(u):=\nabla f(\bar x)u+\nabla^2 f(\bar x)(d,d)+\delta_{T_C^2(\bar x;d)}(u)$ for all $u\in \mathbb{R}^n$. Then, for any $k\in \mathbb{N}$, there exists $u_k\in T_C^2(\bar x;d)$ such that $\varphi(u_k)<\alpha+1/k$.
By Ekeland's variational principle {(\cite[Theorem 3.22]{shap})}, there exists $w_k\in B_{1+\|u_k\|}(u_k)$ such that $\varphi(w_k)\leq \varphi(u_k)<\alpha+1/k$
and that $\varphi{(\cdot)}+\frac{1}{k(1+\|u_k\|)}\|\cdot-w_k\|$ attains its minimum value at $w_k$.
Then, $w_k\in T_C^2(\bar x;d)$ and $\varphi(w_k)=\nabla f(\bar x)w_k+\nabla^2 f(\bar x)(d,d)\rightarrow \alpha$.
By the first-order {optimality} condition, we have
$$0\in \partial \left(\varphi+\frac{1}{k(1+\|u_k\|)}\|\cdot-w_k\|\right)(w_k)\subseteq \nabla f(\bar x)+N_{T_C^2(\bar x;d)}(w_k)+ \frac{1}{k(1+\|u_k\|)}B_{\mathbb{R}^n}.$$
Then, there exists $b_k\in B_{\mathbb{R}^n}$ such that $-\nabla f(\bar x)+b_k/k(1+\|u_k\|)\in N_{T_C^2(\bar x;d)}(w_k)$.
By MSCQ and Proposition \ref{sub}, there exists $\lambda_k\in N_{T_K^{2}(g(\bar x);\nabla g(\bar x)d)}(\nabla g(\bar x)w_k+\nabla^2 g(\bar x)(d,d))$ such that $-\nabla f(\bar x)+b_k/k(1+\|u_k\|)=\nabla g(\bar x)^T\lambda_k$ and
$\|\lambda_k\|\leq \kappa\|\nabla g(\bar x)^T\lambda_k\|=\kappa\|-\nabla f(\bar x)+b_k/k(1+\|u_k\|)\|$. Then, $\{\lambda_k\}$ is bounded.
Without loss of generality, we may assume that $\lambda_k\rightarrow \lambda$. Passing to the limit, we obtain that $-\nabla f(\bar x)=\nabla g(\bar x)^T\lambda$
and $\|\lambda\|\leq \kappa\|\nabla f(\bar x)\|$, i.e., \eqref{oop} holds. Therefore, { The conclusion (ii) of }Theorem \ref{nss} holds true. We also have $$|\langle\nabla g(\bar x)^T\lambda_k+\nabla f(\bar x),w_k\rangle|\leq\frac{\|w_k\|}{k(1+\|u_k\|)}\leq\frac{\|w_k-u_k\|+\|u_k\|}{k(1+\|u_k\|)}\leq\frac{2}{k}\rightarrow 0.$$
In particular, it is easy to observe from Lemma \ref{xxzz} that
$\lambda_k\in N_{T_K^{2}(g(\bar x);\nabla g(\bar x)d)}(\nabla g(\bar x)w_k+\nabla g(\bar x)(d,d))\subseteq N_K(g(\bar x);\nabla g(\bar x)d)$, which indicates that $\lambda\in N_K(g(\bar x);\nabla g(\bar x)d)$. Therefore, $\lambda\in \Lambda(\bar x;d)$.
\endproof

Recall that in the recent work \cite{ye22}, the authors introduced the concept of lower generalized support function as follows:
\begin{definition}\cite[Definition 7]{ye22}\label{gsf}
Let $S\subset \mathbb{R}^n$ be a nonempty closed set. For every subset $A\subset \mathbb{R}^n$, the lower generalized support function to $S$ with respect to $A$ is defined as the mapping {$\hat \sigma_{S,A} : \mathbb{R}^n\rightarrow \mathbb{R}\cup\{\pm\infty\}$} by
$$\hat \sigma_{S,A}(\lambda):=\liminf_{\lambda'\rightarrow\lambda}\inf\limits_{u}\{\langle\lambda',u\rangle \mid u\in N_S^{-1}(\lambda')\cap A\}.$$
If $S=\emptyset$, then we define $\hat \sigma_{S,A}(\lambda)=-\infty$ for all $\lambda$. When $A=\mathbb{R}^n$, we use $\hat \sigma_{S}$ in place of $\hat \sigma_{S,\mathbb{R}^n}$.
\end{definition}

Note that by the definition, we have ${\hat\sigma}_S \leq
\hat \sigma_{S,A}$ for every subset $A\subseteq \mathbb{R}^n$ and $\hat \sigma_{S,B} \leq  \hat \sigma_{S,A}$ whenever $A\subseteq B \subseteq \mathbb{R}^n$. We also observe that $\hat \sigma_{S,A}(\lambda)=+\infty$ {if} for all $\lambda'$ close to $\lambda$, $N_S^{-1}(\lambda')\cap A=\emptyset$.

The additional properties for the multipliers in Theorem \ref{lmby} allow us to derive the following relationships.
\begin{proposition}\label{remark4.6}Let $d\in T_C(\bar x)$ and $\bar x$ be a local optimal solution of problem \eqref{op} in direction $d$. Suppose that $\nabla f(\bar x)d=0$ and MSCQ holds {at $\bar x$} in direction $d$ for the constraint system $g(x)\in K$. Then the following relationships hold:
\begin{enumerate}[{\rm (i)}]
\item 
There exists $\lambda\in \Lambda(\bar x;d)$ such that for any $A \supseteq \Theta$, we have
$$\hat\sigma_{T^{''}_K (g(\bar x);\nabla g(\bar x)d)}(\lambda) \leq\hat\sigma_{T^{''}_K (g(\bar x);\nabla g(\bar x)d),A}(\lambda) \leq \sigma_{\Theta}(\lambda).$$
\item If $T^2_K (g(\bar x);\nabla g(\bar x)d)\not=\emptyset$, then there exists $\lambda\in \Lambda(\bar x;d)$ such that,
for any $B \supseteq \Omega$, we have
$$\hat{\sigma}_{T^2_K (g(\bar x);\nabla g(\bar x)d)}(\lambda)\leq\hat\sigma_{T^2_K (g(\bar x);\nabla g(\bar x)d),B}(\lambda) \leq \sigma_{\Omega}(\lambda).$$
\end{enumerate}
\end{proposition}

\beginproof 
(i)By Theorem \ref{lmby}(i), there exists $\lambda\in \Lambda(\bar x;d)$ such that $\lambda\in N_{T_K^{''}(g(\bar x);\nabla g(\bar x)d)}{(0)}$. Therefore, for any $A\supseteq \Theta$, we have 
\begin{eqnarray} \label{inpp}
\aligned
\hat\sigma_{T^{''}_K (g(\bar x);\nabla g(\bar x)d)}(\lambda)&\leq\hat\sigma_{T^{''}_K (g(\bar x);\nabla g(\bar x)d),A}(\lambda)\\
& := \liminf_{\lambda'\rightarrow\lambda}\inf\limits_{u\in A}\{\langle\lambda',u\rangle \mid \lambda' \in N_{T^{''}_K (g(\bar x);\nabla g(\bar x)d)}(u)\}\\
&\leq \langle\lambda,0\rangle 
\leq \sigma_{\Theta}(\lambda).
\endaligned
\end{eqnarray}

(ii)It follows from Theorem \ref{lmby} (ii) that, there exist  
$w_k\in T_C^2(\bar x;d)$ and $\lambda_k\in N_{T_K^{2}(g(\bar x);\nabla g(\bar x)d)}(\nabla g(\bar x)w_k+\nabla^2 g(\bar x)(d,d))$ 
such that $\lambda_k$ converges to some $\lambda\in \Lambda(\bar x;d)$, $\langle\nabla g(\bar x)^T\lambda_k+\nabla f(\bar x),w_k\rangle\rightarrow 0$.
 Since $\lambda\in \Lambda(\bar x;d)$, it follows that 
 $$\langle\lambda_k-\lambda,\nabla g(\bar x)w_k\rangle=\langle\nabla g(\bar x)^T\lambda_k+\nabla f(\bar x),w_k\rangle\rightarrow 0.$$
Therefore for any $B\supseteq \Omega$, we have $\nabla g(\bar x)w_k+\nabla^2 g(\bar x)(d,d)\in \Omega \subseteq B$ and
\begin{eqnarray}\label{ineq}
\aligned
\hat\sigma_{T^2_K (g(\bar x);\nabla g(\bar x)d)}(\lambda)&\leq\hat\sigma_{T^2_K (g(\bar x);\nabla g(\bar x)d),B}(\lambda)\\&\leq \liminf_{k\rightarrow\infty}\langle\lambda_k,\nabla g(\bar x)w_k+\nabla^2 g(\bar x)(d,d)\rangle 
\\&\leq \liminf_{k\rightarrow\infty}(\sigma_{\Omega}(\lambda)+\langle\lambda_k-\lambda,\nabla g(\bar x)w_k+\nabla^2 g(\bar x)(d,d)\rangle)
\\&=\sigma_{\Omega}(\lambda),
\endaligned
\end{eqnarray}
which completes our proof.
\endproof

Combining Theorem \ref{lmby} with Proposition \ref{remark4.6} we have the following result. Note that the obtained second-order necessary optimality condition utilizing the lower generalized support funtion is  strongest when the set $A$ satisfying 
$A\supseteq \Omega$ is taken as $A=\Omega$ and the set $B$ satisfying 
$B\supseteq \Theta$ is taken as $B=\Theta$. However, even in this case, it is still weaker than the necessary condition  in Theorem \ref{lmby} by virtue of Proposition \ref{remark4.6}.
\begin{theorem}\label{cclmby}
Let $d\in T_C(\bar x)$ and $\bar x$ be a local optimal solution of problem \eqref{op} in direction $d$. Suppose that $\nabla f(\bar x)d=0$ and MSCQ holds {at $\bar x$} in direction $d$ for the constraint system $g(x)\in K$. Then the following conditions hold: 
\begin{enumerate}[{\rm (i)}]
\item There exists $\lambda\in \Lambda(\bar x;d)$ such that { for any set $A\supseteq \Theta$, we have $\hat\sigma_{T^{''}_K (g(\bar x);\nabla g(\bar x)d),A}(\lambda)\leq  0$. In particular 
$$\hat\sigma_{T^{''}_K (g(\bar x);\nabla g(\bar x)d)}(\lambda)\leq  0;$$}

\item If $T^2_K (g(\bar x);\nabla g(\bar x)d)\not=\emptyset$, then there exists $\lambda\in \Lambda(\bar x;d)$ such that, for any set $B\supseteq \Omega$, we have
\begin{equation}\label{aayjj}
\nabla_{xx}^2L(\bar x,\lambda)(d,d)-\hat\sigma_{T^2_K (g(\bar x);\nabla g(\bar x)d),B}(\lambda)\geq 0.
\end{equation}
In particular, 
\begin{equation}\label{aayjjnew}
\nabla_{xx}^2L(\bar x,\lambda)(d,d)-\hat\sigma_{T^2_K (g(\bar x);\nabla g(\bar x)d)}(\lambda)\geq  0.
\end{equation}
\end{enumerate}
\end{theorem}

\begin{remark}\label{fg3}
In addition to the blanket assumption that $\bar x$ is a local optimal solution of problem \eqref{op} only in direction $d$, Theorem \ref{cclmby} can be viewed as an improvement of \cite[Theorem 2]{ye22} in that  Theorem \ref{cclmby} (i) provides supplementary information. Note that Theorem \ref{cclmby} (i) holds regardless whether or not $T^2_K (g(\bar x);\nabla g(\bar x)d)$ is empty or not. In the case where $T^2_K (g(\bar x);\nabla g(\bar x)d)$ is empty, (\ref{aayjjnew}) holds automatically and hence does not provide any information while  Theorem \ref{cclmby} (i) provides some information.
\end{remark}

Recall that the directional Clarke multiplier set:
$$\Lambda^c(\bar x;d):=\{\lambda\in \mathbb{R}^n \mid \lambda\in N^c_K(g(\bar x);\nabla g(\bar x)d)\,\,{\rm and}\,\,\nabla_x L(\bar x,\lambda)=0\}$$
is closed and convex, and generally it is larger than the directional Mordukhovich multiplier set $\Lambda(\bar x;d)$. In the sequel, under the directional Robinson's constraint qualification (DirRCQ):
\begin{equation}\label{23}
\nabla g(\bar x)^T\lambda=0,\lambda\in N^c_K(g(\bar x);\nabla g(\bar x)d)\Rightarrow \lambda=0
\end{equation}
which is stronger than the directional {first-order sufficient condition for metric subregularity (FOSCMS)(\cite[Theorem 1]{gferkl16}) which means that (\ref{23}) holds with $N^c_K(g(\bar x);\nabla g(\bar x)d)$ replaced by the directional limiting normal cone  $N_K(g(\bar x);\nabla g(\bar x)d)$}, we integrate the approach of convex duality theory employed in the recent work of Gfrerer, Ye and Zhou \cite{ye22} to derive the anticipated new second-order necessary optimality condition for problem \eqref{op} in terms of directional Clarke multipliers.

\begin{proposition}\label{pcc}
Let $d\in T_C(\bar x)$ and $\bar x$ be a local optimal solution of problem \eqref{op} in direction $d$ with $\nabla f(\bar x)d=0$. {Suppose that  DirRCQ} holds in direction $d$.
Then, the following three statements are equivalent:
\begin{enumerate}[{\rm (i)}]
\item The primal second-order necessary condition $$\nabla f(\bar x)v \geq 0,\quad\forall v\in T_C^{''}(\bar x,d)$$
of Proposition \ref{op2} holds; 

\item For every $u\in T^{''}_K (g(\bar x);\nabla g(\bar x)d)$, there exists $\lambda_u\in \Lambda^c(\bar x;d)$ such that
$\langle\lambda_u,u\rangle\leq 0$;
\item For every nonempty convex subset $K(d)\subseteq T^{''}_K (g(\bar x);\nabla g(\bar x)d)$, there exists $\lambda\in \Lambda^c(\bar x;d)$ such that
$\sigma_{K(d)}(\lambda)\leq 0$.
\end{enumerate}
\end{proposition}

\beginproof
Since DirRCQ implies that  MSCQ holds in direction $d$ for the constraint system $g(x)\in K$, it follows from Proposition \ref{qc} that
$$T_C^{''}(\bar x;d)= \{v\in \mathbb{R}^n \mid \nabla g(\bar x)v\in T_K^{''}(g(\bar x);\nabla g(\bar x)d)\}.$$

First, we show that ``(i)$\Rightarrow$(ii)". Take arbitary $u\in T_K^{''}(g(\bar x);\nabla g(\bar x)d)$. By Proposition \ref{ppzz}, we have
$u+\hat T_K(g(\bar x);\nabla g(\bar x)d)\subseteq T_K^{''}(g(\bar x);\nabla g(\bar x)d)$.
Assumption (i) ensures that the following conic linear program
\begin{equation}\label{clp}
\min_{v} \nabla f(\bar x)v\,\,\, {\rm s.t.}\,\,\,\nabla g(\bar x)v\in u+\hat T_K(g(\bar x);\nabla g(\bar x)d)
\end{equation}
has nonnegative optimal value. The dual program of the conic linear program \eqref{clp} {(c.f. \cite[Section 2.5.6]{shap})} is
\begin{equation}\label{dclp}
\max_{\lambda\in (\hat T_K(g(\bar x);\nabla g(\bar x)d))^\circ} -\lambda^Tu\,\,\quad {\rm s.t.}\quad\,\,\nabla g(\bar x)^T\lambda+\nabla f(\bar x)=0.
\end{equation}
Since $(\hat T_K(g(\bar x);\nabla g(\bar x)d))^\circ=N^c_K(g(\bar x);\nabla g(\bar x)d)$,  the preceding dual problem can be equivalently written as
$$\max_{\lambda\in \Lambda^c(\bar x;d)} -\lambda^Tu.$$
By \cite[Lemma 6]{ye22}, the DirRCQ implies that
$$0\in {\rm int}(\nabla g(\bar x)\mathbb{R}^n-u-\hat T_K(g(\bar x);\nabla g(\bar x)d)),$$
which indicates that there is no dual between problems \eqref{clp} and \eqref{dclp} according to \cite[Theorem 2.187]{shap}. Hence, the dual program has an optimal solution $\lambda_u$ such that
$$-\langle\lambda_u,u\rangle=\max_{\lambda\in \Lambda^c(\bar x;d)} -\lambda^Tu\geq 0,$$
which establishes the validity of (ii).

Conversely, take $v\in T_C^{''}(\bar x;d)$ and set $u:=\nabla g(\bar x)v$. Then $u\in T_K^{''}(g(\bar x);\nabla g(\bar x)d)$, and it follows from (ii) that there exists $\lambda_u\in \Lambda^c(\bar x;d)$ satisfying 
$\langle\lambda_u,u\rangle\leq 0$. Therefore
$$\nabla f(\bar x)v=-\langle g(\bar x)^T\lambda_u,v\rangle=-\langle\lambda_u,u\rangle\geq 0.$$
which shows that ``(ii)$\Rightarrow$(i)".

To show that ``(ii)$\Leftrightarrow$(iii)", it suffices to establish that ``(ii)$\Rightarrow$(iii)" since the converse implication holds trivially. To this end, pick any nonempty convex subset
${K(d)}\subseteq T^{''}_K (g(\bar x);\nabla g(\bar x)d)$. Without loss of generality, we may assume that { $K(d)$} is closed since $\sigma_{K(d)}=\sigma_{{\rm cl}(K(d))}$. For each $u\in {K(d)}$, assumption (ii) ensures that $\Lambda^c(\bar x;d)$ is nonempty, and then we conclude from \cite[Lemma 6]{ye22} that $\Lambda^c(\bar x;d)$ is
compact. Therefore, it follows from {\cite[Corollary 37.3.2]{rock0}} that
$$\inf_{u\in {K(d)}}\sup_{\lambda\in \Lambda^c(\bar x;d)} -\lambda^Tu=\sup_{\lambda\in \Lambda^c(\bar x;d)}\inf_{u\in {K(d)}} -\lambda^Tu
=\sup_{\lambda\in \Lambda^c(\bar x;d)}-\sigma_{{K(d)}}(\lambda).$$
Recall that the left-hand side of the above equality has nonnegative value according to (ii), and the supremum of the right-hand side of the above equality is attained at some $\lambda$, we conclude that (iii) holds true.
\endproof

Utilizing Proposition \ref{pcc}, Proposition \ref{op2} and \cite[Corollary 4]{ye22}, we easily obtain the following second-order necessary optimality condition which is an improvement of the classical second-order necessary optimality condition (\ref{spc}) in that $K$ can be nonconvex and $T_K^2(g(\bar{x});\nabla g(\bar x)d)$ can be empty.

\begin{corollary}\label{fffd}
Let $d\in T_C(\bar x)$ and $\bar x$ be a local optimal solution of problem \eqref{op} in direction $d$ with $\nabla f(\bar x)d=0$. Suppose that the DirRCQ holds in direction $d$. Then the following conditions hold:
\begin{enumerate}[{\rm (i)}]
\item 
For every nonempty convex {subset} ${K(d)} \subseteq T^{''}_K (g(\bar x);\nabla g(\bar x)d)$, there exists $\lambda\in \Lambda^c(\bar x;d)$ such that
$\sigma_{{K(d)}}(\lambda)\leq 0$.
\item If $T^2_K (g(\bar x);\nabla g(\bar x)d)\not=\emptyset$, then for every nonempty convex {subset} ${K(d)}\subset T^{2}_K (g(\bar x);\nabla g(\bar x)d)$, there exists $\lambda\in \Lambda^c(\bar x;d)$ such that $\nabla_{xx}^2L(\bar x,\lambda)(d,d)-\sigma_{K(d)}(\lambda)\geq0$.
\end{enumerate}
\end{corollary}

Recall that according to \cite[Lemma 7]{ye22}, the directional nondegeneracy condition:
\begin{equation}\label{31}
\nabla g(\bar x)^T\lambda=0,\lambda\in {\rm span} (N_K(g(\bar x);\nabla g(\bar x)d))\Rightarrow \lambda=0
\end{equation}
ensures that $\Lambda(\bar x;d)=\Lambda^c(\bar x;d)$ is a singleton. Hence by Proposition \ref{op2}, Proposition \ref{pcc}(ii) and \cite[Proposition 8]{ye22}, we finally obtain the following sharper necessary second-order condition.

\begin{corollary}\label{fffdss}
Let $d\in T_C(\bar x)$ and $\bar x$ be a local optimal solution of problem \eqref{op} in direction $d$ with $\nabla f(\bar x)d=0$. Suppose that the directional nondegeneracy condition \eqref{31} holds. Then there exists a unique multiplier $\lambda$ such that $\Lambda(\bar x;d) = \{\lambda\}$ which satisfies the following:
\begin{enumerate}[{\rm (i)}]
\item 
$\sigma_{T^{''}_K (g(\bar x);\nabla g(\bar x)d)}(\lambda) {\leq} 0$;
\item If $T^2_K (g(\bar x);\nabla g(\bar x)d)\not=\emptyset$, then $\nabla_{xx}^2L(\bar x,\lambda)(d,d)-\sigma_{T^{2}_K (g(\bar x);\nabla g(\bar x)d)}(\lambda)\geq0$.
\end{enumerate}
\end{corollary}

The following example provides an illustration for the aforementioned results. 
\begin{example}\label{ssa}
Consider the simple constrained optimization problem
\begin{equation*}
\aligned
&\min\ \  f(x_1,x_2):=(x_2+1)^2\\
&{\rm s.t.}\,\,g(x_1,x_2):=(x_1,x_2)\in K,
\endaligned
\end{equation*}
where $K:=\{(y_1,y_2)\in \mathbb{R}^2:y_2=|y_1|^{\frac{3}{2}}\}$. Let $\bar x=(0,0)\in \mathbb{R}^2$ and $d=(1,0)\in \mathbb{R}^2$. Direct calculation yields that $N_K(g(\bar x);\nabla g(\bar x)d)=\{0\}\times \mathbb{R}$, $T^2_K (g(\bar x);\nabla g(\bar x)d)=\emptyset$ and $T^{''}_K (g(\bar x);\nabla g(\bar x)d)=\mathbb{R}\times \mathbb{R}_+$. Then, it is easy to observe that the directional nondegeneracy condition \eqref{31} holds. Hence assertions (i) and (ii) of Corollary \ref{fffdss} hold true. Indeed, $\Lambda(\bar x;d) = \{\lambda\}=\{(0,-2)\}$ such that
$\sigma_{T^{''}_K (g(\bar x);\nabla g(\bar x)d)}(\lambda)= 0$.
\end{example}

Now we are ready to develop second-order sufficient conditions for problem \eqref{op} to obtain a local optimal solution in direction $d$ and show that there is ``no-gap"  in comparison with the necessary conditions established in Theorem \ref{nss}. To this end, we need the following fact which was established in \cite{xde}.

\begin{lemma}[\cite{xde}, Lemma 3.4]\label{li}
Let $\bar x\in S\subseteq \mathbb{R}^n$. If the sequence $x_k\in S\backslash\{\bar x\}$
converges to $\bar x$ such
that $t^{-1}_k(x_k-\bar x)$ converges to some {nonvanishing} vector $d\in T_S(\bar x)$, where $t_k=\|x_k-\bar x\|$.
Then, either $\frac{x_k-\bar x-t_kd}{\frac{1}{2}t_k^2}$ converges to some vector $w\in T_S^2(\bar x;d)\cap\{d\}^\bot$, or
there exists a sequence $r_k\downarrow 0$ such that $t_k/r_k\rightarrow 0$ and $\frac{x_k-\bar x-t_kd}{\frac{1}{2}t_kr_k}$ converges to some vector $w\in T_S^{''}(\bar x;d)\cap \left (\{d\}^\bot\backslash\{0\}\right )$.
\end{lemma}

\begin{theorem}[Second-Order Sufficient Optimality Condition]\label{sff}
Let $d\in T_C(\bar x)\backslash\{0\}$ and $\bar x$  be a feasible solution of problem \eqref{op}.
Furthermore, assume that $\nabla f(\bar x)d=0$ and there is $\lambda\in \mathbb{R}^m$ with $\nabla_x L(\bar x,\lambda)=0$, such that
\begin{enumerate}[{\rm (i)}]
\item $\langle\lambda,v\rangle< 0,\quad\forall v\in \nabla g(\bar x)\left(T_C^{''}(\bar x;d)\cap \left (\{d\}^\bot\backslash\{0\}\right )\right)$;
\item $\nabla_{xx}^2L(\bar x,\lambda)(d,d)-\sigma_{\nabla g(\bar x)(T_C^2(\bar x;d)\cap\{d\}^\bot)+\nabla^2 g(\bar x)(d,d)}(\lambda)>0$.
\end{enumerate}
Then, the second-order growth condition holds at $\bar x$ in direction $d$, that is, there exist constants $\kappa,\rho,\delta>0$ such that
\begin{equation}\label{erjj}
f(x)\geq f(\bar x)+\kappa\|x-\bar x\|^2,\,\,\,\forall x\in C\cap (\bar x+V_{\rho,\delta}(d)),
\end{equation}
which implies that   $\bar x$ is a local minimizer of \eqref{op} in direction $d$. In particular,  suppose that $\nabla f(\bar x)d\geq 0$ for all $d\in T_C(\bar{x})$ and there exists $\lambda\in \mathbb{R}^m$ with $\nabla_x L(\bar x,\lambda)=0$ such that conditions (i) and (ii) in Theorem \ref{sff} hold  for every  nonzero critical direction $d\in C(\bar x)\setminus \{0\}$,  then  $\bar x$ is a local optimal solution of problem \eqref{op}.
\end{theorem}

\beginproof
We argue by contradiction. Assume to the contrary that the second-order growth condition in direction $d$ does not hold at $\bar x$. Then, there {exists a} sequence
$\{x_k\}\subseteq \bar x+V_{1/k,1/k}(d)$ such that $g(x_k)\in K$ and
\begin{equation}\label{fc}
f(x_k)<f(\bar x)+\frac{1}{k}\|x_k-\bar x\|^2.
\end{equation}
According to the assumptions, let $\lambda\in \mathbb{R}^m$ with $\nabla_x L(\bar x,\lambda)=0$ be such that (i) and (ii) hold. Then, we have
$\nabla f(\bar x)+\nabla g(\bar x)^T\lambda=0$.
Let $t_k=\|x_k-\bar x\|$ and $d_k=\frac{x_k-\bar x}{t_k}$. Since $x_k\in \bar x+V_{1/k,1/k}(d)$, one has $t_k\downarrow 0$ and $d_k\rightarrow \frac{d}{\|d\|}\in T_C(\bar x)$.
By Lemma \ref{li}, we have either

(a) $w_k:=\frac{x_k-\bar x-t_k\frac{d}{\|d\|}}{\frac{1}{2}t_k^2}$ converges to some vector $w\in T_C^2(\bar x;\frac{d}{\|d\|})\cap\{d\}^\bot$ \\or

(b) there exists a sequence $r_k\downarrow 0$ such that $t_k/r_k\rightarrow 0$ and $\tilde{w}_k:=\frac{x_k-\bar x-t_k\frac{d}{\|d\|}}{\frac{1}{2}t_kr_k}$ converges to some vector $\tilde{w}\in T_C^{''}(\bar x;\frac{d}{\|d\|})\cap \left (\{d\}^\bot\backslash\{0\}\right )$.

If condition (a) holds, in this case, we have $x_k=\bar x+t_k\frac{d}{\|d\|}+\frac{1}{2}t_k^2w_k$.
Note that $\nabla f(\bar x)d=0$, it follows from \eqref{fc} that
$$\frac{1}{k}t_k^2 > f(x_k)-f(\bar x)= t_k\nabla f(\bar x)\frac{d}{\|d\|}+\frac{1}{2}t_k^2\left(\nabla f(\bar x)w_k+\nabla^2 f(\bar x)\left(\frac{d}{\|d\|},\frac{d}{\|d\|}\right)\right)+o(t_k^2).$$
Dividing both sides of above inequality by $\frac{1}{2}t_k^2$ and passing to the limit, we have
\begin{equation}\label{ww}
\nabla f(\bar x)w+\nabla^2 f(\bar x)\left(\frac{d}{\|d\|},\frac{d}{\|d\|}\right)\leq 0.
\end{equation}
For convenience, let $\Lambda:=\nabla g(\bar x)(T_C^2(\bar x;d)\cap\{d\}^\bot)+\nabla^2 g(\bar x)(d,d)$. Since $w\in T_C^2(\bar x;\frac{d}{\|d\|})\cap\{d\}^\bot$, by \eqref{qctt}, we have
\begin{eqnarray*}
\nabla g(\bar x)w+\nabla^2 g(\bar x)\left(\frac{d}{\|d\|},\frac{d}{\|d\|}\right)&\in& \nabla g(\bar x)\left(T_C^2\left(\bar x;\frac{d}{\|d\|}\right)\cap\{d\}^\bot\right)+\frac{1}{\|d\|^2}\nabla^2 g(\bar x)(d,d)\\&=&\frac{1}{\|d\|^2}\Lambda.
\end{eqnarray*}
And then
\begin{equation}\label{ww1}
\left\langle\lambda,\nabla g(\bar x)w+\nabla^2 g(\bar x)\left(\frac{d}{\|d\|},\frac{d}{\|d\|}\right)\right\rangle
\leq \sigma_{\frac{1}{\|d\|^2}\Lambda}(\lambda)=\frac{1}{\|d\|^2}\sigma_{\Lambda}(\lambda).
\end{equation}
Adding \eqref{ww} and \eqref{ww1}, we obtain that
\begin{equation*}
\aligned
&\nabla f(\bar x)w+\langle\nabla g(\bar x)^T\lambda,w\rangle+\nabla^2 f(\bar x)\left(\frac{d}{\|d\|},\frac{d}{\|d\|}\right)+
\left\langle\lambda,\nabla^2 g(\bar x)\left(\frac{d}{\|d\|},\frac{d}{\|d\|}\right)\right\rangle\\
&= 0+\nabla^2_{xx} L(\bar x,\lambda)\left(\frac{d}{\|d\|},\frac{d}{\|d\|}\right)=\frac{1}{\|d\|^2}\nabla^2_{xx} L(\bar x,\lambda)\left(d,d\right) \leq \frac{1}{\|d\|^2}\sigma_{\Lambda}(\lambda),
\endaligned
\end{equation*}
which is a contradiction to assumption (ii).

If condition (b) holds, then we have $x_k=\bar x+t_k\frac{d}{\|d\|}+\frac{1}{2}t_kr_k\tilde{w}_k$, {and}  ${\tilde w} \in T''_C(\bar x, d)$ according to \eqref{qctt}.
It follows from \eqref{fc} that
\begin{eqnarray*}
\frac{1}{k}t_k^2  &>& f(x_k)-f(\bar x)\\
&=& t_k\nabla f(\bar x)\frac{d}{\|d\|}+\frac{1}{2}t_kr_k\nabla f(\bar x)\tilde{w}_k+\frac{1}{2}t_k^2\nabla^2 f(\bar x)\left(\frac{d}{\|d\|},\frac{d}{\|d\|}\right)+o(t_k^2).
\end{eqnarray*}
Since $t_k/r_k\rightarrow 0$, dividing both sides of above inequality by $\frac{1}{2}t_kr_k$ and passing to the limit, we obtain that
$-\langle \lambda, \nabla g(\bar x)\tilde{w}\rangle=\nabla f(\bar x)\tilde{w}\leq 0,$ which is a contradiction to (i). Hence, we conclude that the second-order growth condition holds at $\bar x$ in direction $d$.
\endproof

By the relationships (\ref{strin})-(\ref{strintheta}), conditions in Theorem \ref{sff} (i) and (ii) can be replaced by some stronger but  easier to verify conditions as in the following corollary.
\begin{corollary}\label{cor4.15} {Let $d\in T_C(\bar x)\setminus\{0\}$ and $\bar x$  be a feasible solution of problem \eqref{op}.
Furthermore, assume that $\nabla f(\bar x)d=0$ and there is $\lambda\in \mathbb{R}^m$ with $\nabla_x L(\bar x,\lambda)=0$, such that
\begin{enumerate}[{\rm (i)}]
\item $\langle\lambda,v\rangle< 0,\quad\forall v\in  T''_K (g(\bar x);\nabla g({\bar x})d)  \cap\nabla g(\bar x) \left (\{d\}^\bot\setminus\{0\}\right )$;
\item $ \nabla_{xx}^2L(\bar x,\lambda)(d,d)-\sigma_{T^2_K (g(\bar x);\nabla g({\bar x})d)}(\lambda)>0$.
\end{enumerate}
Then, the second-order growth condition holds at $\bar x$ in direction $d$. In particular suppose that  $\nabla f(\bar x)d\geq 0$ for all $d \in T_C(\bar x)$   and  and there exists $\lambda\in \mathbb{R}^m$ with $\nabla_x L(\bar x,\lambda)=0$ such that conditions (i) and (ii) above hold for every nonzero critical direction $d \in C(\bar x)\setminus \{0\}$, then the second-order growth condition holds at $\bar x$.}
\end{corollary}

We can observe from Example \ref{fg3a} that, the value of $\sigma_{T_K^2(g(\bar x);\nabla g(\bar x)d)}(\lambda)$ can be strictly greater than the value of $\sigma_{\nabla g(\bar x)(T_C^2(\bar x;d)\cap\{d\}^\bot)+\nabla^2 
g(\bar x)(d,d)}(\lambda)$. Hence, the sufficient condition in Theorem \ref{sff}(ii)  is strictly weaker than the one in Corollary \ref{cor4.15}(ii).
As commented in the introduction, in \cite[Theorem 3.86]{shap} and \cite[Theorem 4]{ye22}, under the assumption of outer second-order regularity on $K$ and the convexity of $T_K^2(g(\bar x);\nabla g(\bar x)d)$, the authors provide a second-order sufficient condition in the form of (\ref{spcsuffcl}).
Besides, Theorem \ref{sff} and Corollary \ref{cor4.15} also give the exact relationship between the second-order sufficient condition and the corresponding direction considered in the second-order growth condition.

We provide the following example to illustrate our theoretical results in  Corollary \ref{cor4.15}.

\begin{example}\label{ccff}
Consider the simple constrained optimization problem
\begin{equation*}
\aligned
\min &f(x_1,x_2):=x_1+x_1^2+x_2^2\\
&{\rm s.t.}\,\,g(x_1,x_2):=(x_1+x_2^2,2x_2)\in K:=\{(y_1,y_2)\in \mathbb{R}^2 \mid y_1\geq y_2^2/2\}.
\endaligned
\end{equation*}
Let $\bar x=(0,0)\in \mathbb{R}^2$ and $d=(0,1)\in \mathbb{R}^2$. 
 Clearly, the feasible region $C:=g^{-1}(K)=\{(x_1,x_2)\in \mathbb{R}^2 \mid x_1\geq x_2^2\}$ and $d\in T_C(\bar x)={\mathbb{R}_+\times \mathbb{R}}$.
By direct calculation, we have
$\nabla f(\bar x)=(1,0),\nabla^2 f(\bar x)=\left(
\begin{array}{cc}
 2 & 0\\
 0 & 2\\
\end{array}
\right)
$
and $\nabla g(\bar x)=\left(
\begin{array}{cc}
 1 & 0\\
 0 & 2\\
\end{array}
\right)
$.
It is easy to calculate that $T_K^2(g(\bar x);\nabla g(\bar x)d)=\{(w_1,w_2)\in \mathbb{R}^2 \mid w_1\geq 4 \}$ and $T_K^{''}(g(\bar x);\nabla g(\bar x)d)\cap \nabla g(\bar x) \left (\{d\}^\bot\setminus\{0\}\right ) =\{(w_1,w_2)\in \mathbb{R}^2 \mid w_1 >0, w_2=0 \}$. Let $\lambda=(-1,0)$. Then $\nabla f(\bar x)d =0$, $\nabla_x L(\bar x,\lambda)=0$ and conditions (i) and (ii) of Corollary \ref{cor4.15} hold. Hence we  conclude that $\bar x$ is a strictly local minimizer in direction $d$ fulfilling the second-order
growth condition.  Moreover, it is easy to verify that for any $d\in T_C(\bar x)$, $\nabla f(\bar x)d \geq 0$. Since conditions (i) and (ii) of Corollary 4.15 hold for all nonzero critical directions we can conclude that $\bar x$ is a strictly local minimizer. 
\end{example}

It is interesting to note that instead of requiring the outer second-order regularity on $K$ which may not be easy to verify, we require  Corollary \ref{cor4.15}(i) to hold.   It is worth to note that, in the case of $T_K^2(g(\bar x);\nabla g(\bar x)d)= \emptyset$, condition (ii) of Corollary \ref{cor4.15} holds automatically, which leaves us to check condition (i) of Corollary \ref{cor4.15} only. For example consider  the problem in  Example \ref{ccff} with $K=\{(y_1,y_2)\in \mathbb{R}^2 \mid y_1\geq 0, y_1^2=y_2^3\}$. Then we have $T_K^2(g(\bar x);\nabla g(\bar x)d)= \emptyset$ and $T_K^{''}(g(\bar x);\nabla g(\bar x)d)=\{(w_1,w_2)\in \mathbb{R}^2 \mid w_1\geq 0\}$.   Condition  Corollary \ref{cor4.15}(i) holds. Then, the second-order growth condition holds at $\bar x$ in direction $d$.


\section{Conclusion}
In this paper, we consider second-order optimality conditions in specific directions for a generalized set-constrained optimization problem. By utilizing both the outer second-order tangent set and the asymptotic second-order tangent cone, we establish new directional necessary and ``no gap" sufficient second-order optimality conditions which do not require convexity or outer second-order regularity of the underlying set or nonemptiness of the outer second-order tangent set. The obtained theoretical results are improvement of the existing ones. In our future work, we may consider  extending our results to problems without second-order differentiability in functions $f$ and/or $g$. We also look forward to apply our results  in numerical algorithms. 

\section*{Acknowledgments}
The authors would like to thank one of the anonymous referees for the   suggestions and comments which inspire us to improve the presentation of the paper.


\end{sloppypar}
\end{document}